\documentclass[11pt]{amsart}
\usepackage[margin=1.15in]{geometry}
\usepackage{amssymb,mathtools,enumitem}
\usepackage[colorlinks,linkcolor=blue,citecolor=blue,urlcolor=blue]{hyperref}

\newtheorem{theorem}{Theorem}[section]
\newtheorem{lemma}[theorem]{Lemma}
\newtheorem{proposition}[theorem]{Proposition}
\newtheorem{corollary}[theorem]{Corollary}
\theoremstyle{definition}
\newtheorem{definition}[theorem]{Definition}
\newtheorem{question}[theorem]{Question}
\theoremstyle{remark}
\newtheorem{remark}[theorem]{Remark}

\newcommand{\rest}{{\upharpoonright}}
\newcommand{\leK}{\le_K}
\newcommand{\w}{\omega}
\newcommand{\cs}{2^{\w}}
\newcommand{\incomp}{\mathrel{\bot}}
\newcommand{\fr}[1]{0.#1}
\setlist[enumerate]{leftmargin=2.2em}

\title[Incompatibility and the existential theory of the $K$-degrees]{Incompatibility and the existential theory\\ of the $K$-degrees}
\author{Yongge Wang}
\address{Department of Cybersecurity, University of North Carolina at Charlotte, Charlotte, NC 28223, USA}
\email{yonwang@charlotte.edu}
\date{}

\begin{document}
\begin{abstract}
Two reals are \emph{compatible} if they have a common upper bound in the $K$-degrees. We study incompatibility and
the lattice-theoretic structure of the $K$-degrees. We give a short proof of the Miller--Yu theorem that the two
halves of a random real are incompatible, and a criterion for incompatibility along computable gauges. If
$R\oplus Y$ is random and $F$ is a computable set with counting function $f$, then the real $R_F$ obtained by
setting the bits of $R$ in $F$ to $0$ is $K$-below every random real when $\sum_n2^{-f(n)}<\infty$, and incompatible
with $Y$ otherwise; so the formula ``$x$ is incompatible with some degree'' does not define randomness. If both
reals are overwritten along sets with at most $\alpha\log n$ elements below $n$, $\alpha<1$, they remain
incompatible: the deficiencies do not add. Every finite configuration of $\leK$ and incompatibility allowed by
two obvious axioms occurs, so the existential theory of the $K$-degrees with incompatibility is decidable. We
describe the lower cones of sparse codes exactly, determine the order among sparse overwritings of a single random
real that is low for $\Omega$, and give an explicit common-source construction of least upper bounds and
greatest lower bounds. We show that the $K$-degrees are not
elementarily equivalent to the Turing degrees or to the vL-degrees, give sufficient conditions for the
undecidability of their first-order theory, and collect open problems.
\end{abstract}
\maketitle

\section{Introduction}\label{sec:intro}

For reals $X,Y\in2^\omega$, Downey, Hirschfeldt and LaForte defined $X\leK Y$ to mean that
$K(X\rest n)\le K(Y\rest n)+O(1)$, where $K$ is prefix-free complexity. The resulting degree structure
$\mathcal D_K$, the \emph{$K$-degrees}, measures how random a real is by the complexity of its initial segments.
Its least element consists of the $K$-trivial reals, random reals are closed upward, and several basic questions
about it have only recently been answered: there is a minimal pair \cite{CM}, there are pairs of degrees with no
meet, non-zero degrees that are the meet of two others, and lattice embeddings that preserve meets and joins
in the whole structure \cite{MS}, there is no maximal degree \cite{Liu}, and every
real has an upper bound of effective Hausdorff dimension one \cite{Zhao}. The structure is not an upper
semilattice: by Miller and Yu \cite{MY}, if $X\oplus Y$ is random then $X$ and $Y$ have no common upper bound.

We call two reals \emph{compatible} if they have a common upper bound, and \emph{incompatible} otherwise. Question~9.8
of Miller and Nies \cite{MN}, whether every real is $K$-below a random real, asks precisely whether every real is
compatible with some random real. This paper studies incompatibility, the joins and meets that do exist in
$\mathcal D_K$, and what this says about the first-order theory of $\mathcal D_K$.

\smallskip
\noindent\textbf{Results.}
\begin{enumerate}[label=(\arabic*)]
\item \emph{Joint dips} (Section~\ref{sec:dip}). Every common upper bound $w$ of $x$ and $y$ forces the pair
$(x\rest N,y\rest N)$ to be compressible by about $\log N$ bits at the lengths $N$ announced by $w$
(Theorem~\ref{thm:basic}). This gives a short proof of the Miller--Yu theorem (Corollary~\ref{cor:MY}) and, along
arbitrary computable gauges, a criterion for incompatibility (Corollary~\ref{cor:crit}).
\item \emph{Overwritten random reals} (Section~\ref{sec:over}). If $R\oplus Y$ is random and $F$ is computable with
counting function $f$, then $R_F$, the real $R$ with the positions in $F$ set to $0$, is $K$-below every random real
when $\sum_n2^{-f(n)}<\infty$, and is incompatible with $Y$ otherwise (Theorem~\ref{thm:dich}). Consequently the
formula ``$x$ is incompatible with some degree'' does not define the random degrees (Corollary~\ref{cor:J}).
\item \emph{Finite configurations} (Section~\ref{sec:masks}). Every finite configuration of $\leK$ and incompatibility
allowed by two obvious axioms occurs, so the existential theory of $(\mathcal D_K,\le,\incomp)$ is decidable
(Corollary~\ref{cor:fin}).
\item \emph{Countable lower cones} (Section~\ref{sec:sparse}). For sparse codes $S_V$ the lower cone is described
exactly, and $\leK$ between sparse codes is a strong reducibility close to $rK$-reducibility
(Proposition~\ref{prop:sparse}, Corollary~\ref{cor:sparse}).
\item \emph{Joins and meets} (Section~\ref{sec:timing}). Merkle and Stephan \cite[Section~9]{MS} embed the
lattice of subsets of $\omega$, and every countable distributive lattice, into these degree structures with meets
and joins preserved, realizing joins and meets as pointwise maxima and minima of complexity functions. We give
an explicit construction of this kind for the Boolean algebra of computable sets modulo finite sets, in which all
reals carry one common random source revealed on different timing schedules (Theorem~\ref{thm:timing}). For a different family of reals, two natural degrees have
upper bounds but the natural candidate is not a least one (Proposition~\ref{prop:nojoin}).
\item \emph{Independent sources} (Section~\ref{sec:source}). If $R\oplus Y$ is random and $F,G$ are computable with
counting functions at most $\alpha\log n$, $\alpha<1$, then $R_F$ and $Y_G$ are incompatible
(Theorem~\ref{thm:two}, Corollary~\ref{cor:two}): the deficiencies do not add. The proof shows that a $K$-upper
bound of a nearly random real carries almost no information that the real lacks
(Proposition~\ref{prop:source}). For a single random source that is low for $\Omega$, overwriting along sparse
computable sets reverses inclusion modulo finite sets (Proposition~\ref{prop:order}).
\item \emph{The first-order theory} (Sections~\ref{sec:fo}--\ref{sec:open}). $\mathcal D_K$ is not elementarily
equivalent to the Turing degrees or to any upper semilattice, nor to the vL-degrees
(Propositions~\ref{prop:noniso} and~\ref{prop:vL}); a minimal $K$-degree, if one exists, lies in a narrow zone of
initial segment complexity (Proposition~\ref{prop:dilute}). We formulate sufficient conditions for undecidability
(Propositions~\ref{prop:nies} and~\ref{prop:target}), prove partial results toward them (Theorem~\ref{thm:Bn}), and
collect open problems, old and new (Section~\ref{sec:open}).
\end{enumerate}

Sections~\ref{sec:other}--\ref{sec:need} contain further remarks on Question~9.8. All reals in our constructions are
explicit; the main tools are elementary counting with prefix-free complexity, symmetry of information, and, in
Section~\ref{sec:source}, a relativized form of Muchnik's theorem on conditional descriptions.

\section{Preliminaries}\label{sec:prelim}

$K$ is prefix-free Kolmogorov complexity with respect to a fixed universal prefix-free machine $U$;
$K(\sigma,\tau)$ is the complexity of a pair and $K(\sigma\mid\tau)$ is conditional complexity.
For a string $\sigma$, $\sigma^*$ is a shortest $U$-program for $\sigma$ (the first one found); conditioning on
$\sigma^*$ is the same as conditioning on $(\sigma,K(\sigma))$.
Natural numbers are identified with strings, and $\log$ is to base $2$.
For $X\in\cs$, $X\rest n$ is the prefix of length $n$, and
$X\leK Y$ means $K(X\rest n)\le K(Y\rest n)+O(1)$ for all $n$.
``Random'' means Martin-L\"of random. By the Levin--Schnorr theorem, $X$ is random iff
$K(X\rest n)\ge n-O(1)$. Hence
\begin{equation}\label{eq:up}
 X\leK Y\ \text{and}\ X\ \text{random}\ \Longrightarrow\ Y\ \text{random}.
\end{equation}
$X\oplus Y$ is the join, $(X\oplus Y)(2n)=X(n)$ and $(X\oplus Y)(2n+1)=Y(n)$, so that
\begin{equation}\label{eq:join}
 K\bigl((X\oplus Y)\rest 2n\bigr)\le K(X\rest n,\,Y\rest n)+O(1).
\end{equation}
We use symmetry of information (G\'acs, Chaitin; see \cite[Ch.~3]{DH}):
\begin{equation}\label{eq:soi}
 K(\sigma,\tau)=K(\sigma)+K(\tau\mid\sigma^*)\pm O(1).
\end{equation}

\begin{lemma}\label{lem:pair}
There is a constant $c$ such that for all $N$ and all strings $\sigma,\tau$ of length $N$,
\[
 K(\sigma,\tau)\le K(\sigma)+K(\tau)-K(N)+c .
\]
\end{lemma}
\begin{proof}
Since $N=|\sigma|$ is computable from $\sigma$, \eqref{eq:soi} gives
$K(\sigma)=K(N,\sigma)\pm O(1)=K(N)+K(\sigma\mid N^*)\pm O(1)$, and likewise for $\tau$.
A prefix-free machine that reads $N^*$, then a shortest program for $\sigma$ given $N^*$, then a shortest program for
$\tau$ given $N^*$, shows
$K(\sigma,\tau)\le K(N)+K(\sigma\mid N^*)+K(\tau\mid N^*)+O(1)=K(\sigma)+K(\tau)-K(N)+O(1)$.
\end{proof}

In words: each of $\sigma,\tau$ pays for the length $N$ separately, but the pair needs to pay for it only once.

\section{Compatibility and Question 9.8}

\begin{definition}
Reals $x,y$ are \emph{compatible} if there is $u$ with $x\leK u$ and $y\leK u$; otherwise they are
\emph{incompatible}, written $x\incomp y$. Let $J(x)$ be the statement $\exists y\,(x\incomp y)$.
Both notions depend only on $K$-degrees, and $J$ is a first-order formula in the language of the partial order.
\end{definition}

Question~9.8 of \cite{MN} asks: \emph{for every $A\in\cs$, is there a random $X$ with $A\leK X$?}
It is still open. Recent related results: Liu \cite{Liu} showed that there is no maximal $K$-degree
(Question~9.7 of \cite{MN}); Barmpalias and Liu \cite{BL} showed that the analogue of Question~9.8 for
$rK$-reducibility fails; Zhao \cite{Zhao} showed that every real has a $K$-upper bound of effective Hausdorff
dimension~$1$.

\begin{proposition}\label{prop:basic}
\begin{enumerate}[label=\textup{(\alph*)}]
\item If $x\leK x'$ and $x\incomp y$, then $x'\incomp y$. So $\{x: J(x)\}$ is closed upward and its complement,
the set of reals compatible with every real, is closed downward.
\item Question~9.8 has a positive answer iff every real is compatible with some random real.
A counterexample is exactly a real that is incompatible with every random real.
\item \textup{(Miller--Yu \cite{MY})} Every random real satisfies $J$.
\item If $J$ held \emph{only} of random reals, Question~9.8 would have a positive answer.
\end{enumerate}
\end{proposition}
\begin{proof}
(a) An upper bound of $x'$ and $y$ is an upper bound of $x$ and $y$.
(b) If $x\leK w$ with $w$ random, then $w$ witnesses that $x$ and $w$ are compatible. Conversely, if $u$ is an
upper bound of $x$ and of a random $w$, then $u$ is random by \eqref{eq:up}, and $x\leK u$.
(c) If $x$ is random and $y$ is random relative to $x$, then $x\oplus y$ is random and $x\incomp y$ by \cite{MY}
(reproved as Corollary~\ref{cor:MY}).
(d) A non-random real would then be compatible with every real, in particular with a random one; apply~(b).
\end{proof}

Theorem~\ref{thm:dich} shows that the hypothesis of (d) is false: $J$ holds of some non-random reals.
Those reals have random upper bounds, so this closes one route to Question~9.8 without bearing on its answer.

\section{Joint dips}\label{sec:dip}

For a real $w$ and $m\in\w$, let $w_m$ be the integer whose binary expansion is $1\,w(0)\,w(1)\cdots w(m-1)$.
Thus $2^m\le w_m<2^{m+1}$, and $w_m$ and $w\rest m$ are computable from each other; in particular
$K(w_m)=K(w\rest m)\pm O(1)$.

\begin{lemma}\label{lem:self}
For all reals $w$ and all $m$: $K(w\rest w_m)\le K(w\rest m)+w_m-m+O(1)$.
\end{lemma}
\begin{proof}
Let $M$ be the machine that on input $p\rho$ runs $U(p)$; if this halts with output $\sigma$, $M$ lets $N$ be
the integer with binary expansion $1\sigma$, reads exactly $N-|\sigma|$ further bits $\rho$, and outputs $\sigma\rho$.
The domain of $M$ is prefix-free. With $p=(w\rest m)^*$ and $\rho=w(m)\cdots w(w_m-1)$ the output is $w\rest w_m$.
\end{proof}

\begin{theorem}[Joint dip]\label{thm:basic}
Suppose $x\leK w$ and $y\leK w$. Then there is $c$ such that for every $m$, writing $N=w_m$,
\[
 K(x\rest N,\;y\rest N)\ \le\ 2N-2m+K(w\rest m)+c\ \le\ 2N-m+K(m)+2c .
\]
\end{theorem}
\begin{proof}
By Lemma~\ref{lem:pair}, the two reductions, Lemma~\ref{lem:self}, and $K(N)=K(w\rest m)\pm O(1)$:
\begin{align*}
K(x\rest N,y\rest N)
 &\le K(x\rest N)+K(y\rest N)-K(N)+O(1)\\
 &\le 2K(w\rest N)-K(w\rest m)+O(1)\\
 &\le 2\bigl(K(w\rest m)+N-m\bigr)-K(w\rest m)+O(1)=2N-2m+K(w\rest m)+O(1).
\end{align*}
The second inequality of the theorem follows from $K(w\rest m)\le m+K(m)+O(1)$.
\end{proof}

\begin{corollary}[Miller--Yu \cite{MY}]\label{cor:MY}
If $x\oplus y$ is random, then $x$ and $y$ have no common upper bound in the $K$-degrees.
\end{corollary}
\begin{proof}
If $w$ were a common upper bound, then by \eqref{eq:join} and Theorem~\ref{thm:basic}, for every $m$ and $N=w_m$,
$2N-O(1)\le K((x\oplus y)\rest 2N)\le K(x\rest N,y\rest N)+O(1)\le 2N-m+K(m)+O(1)$.
So $m-K(m)\le O(1)$ for all $m$, which is false because $K(m)\le 2\log m+O(1)$.
\end{proof}

\begin{remark}
The published proof in \cite{MY} goes through van Lambalgen reducibility. We have not found the argument above
in the literature, but it is short enough that it may well be folklore.
The mechanism is the one in Lemma~\ref{lem:pair}: any real $w$ has, at the length $N=w_m$, complexity $m$ below the
maximum $N+K(N)$, because its own first $m$ bits announce the length. Both $x$ and $y$ must be equally low there,
and two strings of the same length that are both $m$ below the maximum are jointly $m-(K(N)-m)$ below $2N$.
\end{remark}

\begin{remark}\label{rem:quick}
Theorem~\ref{thm:basic} already gives more than Corollary~\ref{cor:MY}: if
$K(x\rest n,y\rest n)\ge 2n-\log n+2\log\log n$ for all large $n$, then $x\incomp y$.
(For $N\in[2^m,2^{m+1})$ this would give $2\log m-1\le K(m)+O(1)$, which fails for large $m$.)
Section~\ref{sec:over} gives the sharp form.
\end{remark}

\section{Gauges}\label{sec:gauge}

The length $w_m$ is one way in which a real announces a length at which it must be simple. Martin-L\"of's
covering argument \cite{ML} provides such lengths for every divergent series.

\begin{definition}\label{def:gauge}
A \emph{gauge} is a computable $g:\w\to\w\cup\{\infty\}$ such that $g(n)\le n$ whenever $g(n)<\infty$, and
$\sum_n 2^{-g(n)}=\infty$ (with $2^{-\infty}=0$). Its \emph{blocks} are the intervals $I_i=[n_i,n_{i+1})$, where
$n_0=0$ and $n_{i+1}$ is the least $n>n_i$ with $\sum_{n_i\le k<n}2^{-g(k)}\ge 1$.
\end{definition}

\begin{lemma}[Cover lemma]\label{lem:cover}
For every gauge $g$ there is a constant $c$ such that for every real $u$ and every $i$ there is $N\in I_i$ with
$g(N)<\infty$ and
\begin{enumerate}[label=\textup{(\alph*)}]
\item $K(u\rest N)\le N+K(i)+c$;
\item $K(u\rest g(N))\le K(N)+K(i)+c$.
\end{enumerate}
\end{lemma}
\begin{proof}
For a string $\sigma$ write $\fr{\sigma}=\sum_{k<|\sigma|}\sigma(k)2^{-k-1}$, a multiple of $2^{-|\sigma|}$ in $[0,1)$.
For $n\le n'$ we have
\begin{equation}\label{eq:mono}
 \fr{(u\rest n)}\ \le\ \fr{(u\rest n')}\ <\ \fr{(u\rest n)}+2^{-n}.
\end{equation}
Fix $i$. Put $\alpha_{n_i}=0$ and $\alpha_{n+1}=\min\{1,\alpha_n+2^{-g(n)}\}$ for $n\in I_i$.
By the choice of $n_{i+1}$, $\alpha_{n_{i+1}}=1$. By induction, $\alpha_{n+1}$ is a multiple of $2^{-n}$
(because $g(n)\le n$), hence $\alpha_n$ is a multiple of $2^{-n}$ for all $n_i\le n\le n_{i+1}$.
For $n\in I_i$ let
\[
 A_n=\{\sigma\in 2^n:\ \alpha_n\le \fr{\sigma}<\alpha_{n+1}\}.
\]

\emph{Claim.} For every real $u$ there is exactly one $N\in I_i$ with $u\rest N\in A_N$.

Existence: let $N$ be the least $n\in I_i$ with $\fr{(u\rest n)}<\alpha_{n+1}$; it exists because
$\alpha_{n_{i+1}}=1$. If $N=n_i$ then $\alpha_N=0\le\fr{(u\rest N)}$. Otherwise, by minimality and \eqref{eq:mono},
$\alpha_N\le\fr{(u\rest (N-1))}\le\fr{(u\rest N)}$.
Uniqueness: suppose $n<n'$ in $I_i$ both qualify. As $\fr{(u\rest n)}$ and $\alpha_{n+1}$ are multiples of $2^{-n}$
and $\fr{(u\rest n)}<\alpha_{n+1}$, we get $\fr{(u\rest n)}\le\alpha_{n+1}-2^{-n}$. By \eqref{eq:mono},
$\fr{(u\rest n')}<\alpha_{n+1}\le\alpha_{n'}$, contradicting $u\rest n'\in A_{n'}$.

Let $N$ be given by the claim. If $g(n)=\infty$ then $A_n=\emptyset$, so $g(N)<\infty$.

(a) Consider the machine that reads a $U$-program for some $i$, computes the numbers $\alpha_n$ for the block $I_i$
as needed, and then reads further bits one at a time; as soon as the string $\rho$ of further bits read so far
has length $n\in I_i$ and lies in $A_n$, it halts with output $\rho$. Its domain is prefix-free. By the claim, on
input $i^*$ followed by $u\rest N$ it halts after reading exactly these bits and outputs $u\rest N$.

(b) Let $k=g(N)\le N$. Since $u\rest N\in A_N$, \eqref{eq:mono} gives
$\alpha_N-2^{-k}<\fr{(u\rest k)}<\alpha_{N+1}$. As $\alpha_{N+1}-\alpha_N\le 2^{-k}$, this open interval has length
at most $2\cdot 2^{-k}$ and contains at most two multiples of $2^{-k}$. So $u\rest k$ is one of at most two strings
that can be computed from $i$ and $N$, and $K(u\rest k)\le K(i,N)+O(1)\le K(N)+K(i)+O(1)$.
\end{proof}

\begin{theorem}[Joint dips along a gauge]\label{thm:gauge}
Let $g$ be a gauge with blocks $I_i$, and suppose $x\leK u$ and $y\leK u$. Then there is $c$ such that for every
$i$ there is $N\in I_i$ with $g(N)<\infty$ and
\[
 K(x\rest N,\;y\rest N)\ \le\ 2N-K(u\rest g(N))+3K(i)+c .
\]
\end{theorem}
\begin{proof}
Take $N$ from Lemma~\ref{lem:cover}. By Lemma~\ref{lem:pair}, the reductions, and (a), (b) of Lemma~\ref{lem:cover},
\[
K(x\rest N,y\rest N)\le 2K(u\rest N)-K(N)+O(1)\le 2\bigl(N+K(i)\bigr)-\bigl(K(u\rest g(N))-K(i)\bigr)+O(1).\qedhere
\]
\end{proof}

\begin{remark}
For $g(n)=\lfloor\log n\rfloor$ ($n\ge1$, and $g(0)=\infty$) the blocks are $I_0=[0,2)$ and $I_i=[2^i,2^{i+1})$ for $i\ge1$, and
the length chosen in block $i$ is $N=u_i$. So Theorem~\ref{thm:basic} is this instance of Theorem~\ref{thm:gauge},
with slightly sharper bookkeeping.
\end{remark}

To use Theorem~\ref{thm:gauge} one needs a gauge that exceeds a given deficiency by a margin larger than $3K(i)$.

\begin{lemma}[Slack]\label{lem:slack}
Let $h:\w\to\w$ be computable with $\sum_n 2^{-h(n)}=\infty$. Then there are a gauge $g$ and a computable,
nondecreasing, unbounded $\ell:\w\to\w$ such that
\begin{enumerate}[label=\textup{(\roman*)}]
\item $g(n)=h(n)+\ell(n)+1$ whenever $g(n)<\infty$;
\item if $I_i$ is the $i$-th block of $g$ and $N\in I_i$, then $i\le\ell(N)$.
\end{enumerate}
\end{lemma}
\begin{proof}
Let $a_n=2^{-h(n)}\in(0,1]$, $S_n=1+\sum_{k<n}a_k$ and $\ell(n)=\lceil\log S_n\rceil$, so $S_n\le 2^{\ell(n)}<2S_n$.
Put $g(n)=h(n)+\ell(n)+1$ if this number is at most $n$, and $g(n)=\infty$ otherwise.

\emph{$g$ is a gauge.} Let $B=\{n:g(n)=\infty\}$. For $n\in B$ we have $h(n)>n-\ell(n)-1$, so
$a_n<2^{\ell(n)+1-n}<4S_n2^{-n}\le 4(n+1)2^{-n}$; hence $\sum_{n\in B}a_n<\infty$.
For $n\notin B$, $2^{-g(n)}=a_n2^{-\ell(n)-1}>a_n/(4S_n)$. Now $\sum_n a_n/S_n=\infty$: for $M<M'$,
$\sum_{n=M}^{M'-1}a_n/S_n\ge (S_{M'}-S_M)/S_{M'}$, which tends to $1$ as $M'\to\infty$ because $S_{M'}\to\infty$;
so every tail of the series has sum at least $1/2$, and the series diverges. Since
$\sum_{n\in B}a_n/S_n\le\sum_{n\in B}a_n<\infty$, we get $\sum_{n\notin B}2^{-g(n)}=\infty$.

\emph{Proof of \textup{(ii)}.} Each block has $\sum_{k\in I_j}2^{-g(k)}\ge1$, so
$i\le\sum_{k<n_i}2^{-g(k)}$. For $k\notin B$, $2^{-g(k)}\le a_k/(2S_k)\le a_k/S_{k+1}$, because
$S_{k+1}=S_k+a_k\le 2S_k$. And $a_k/S_{k+1}=1-S_k/S_{k+1}\le\ln(S_{k+1}/S_k)$. Summing,
$i\le\ln S_{n_i}\le\log S_N\le\ell(N)$ for $N\ge n_i$.
\end{proof}

\begin{corollary}[Incompatibility criterion]\label{cor:crit}
Let $h:\w\to\w$ be computable and suppose
\[
 K(x\rest n,\;y\rest n)\ \ge\ 2n-h(n)-O(1)\qquad\text{for all }n.
\]
\begin{enumerate}[label=\textup{(\alph*)}]
\item If $x$ or $y$ is random and $\sum_n2^{-h(n)}=\infty$, then $x\incomp y$.
\item If $h$ is nondecreasing and $\sum_n2^{-2h(n)}=\infty$, then $x\incomp y$ \textup{(}no randomness assumed\textup{)}.
\end{enumerate}
\end{corollary}
\begin{proof}
Suppose $u$ is a common upper bound.

(a) By \eqref{eq:up}, $u$ is random, so $K(u\rest k)\ge k-O(1)$ for all $k$. Let $g,\ell$ be as in
Lemma~\ref{lem:slack} for $h$. By Theorem~\ref{thm:gauge}, for each $i$ there is $N\in I_i$ with
\[
2N-h(N)-O(1)\le K(x\rest N,y\rest N)\le 2N-g(N)+3K(i)+O(1),
\]
that is, $\ell(N)+1=g(N)-h(N)\le 3K(i)+O(1)$. Since $i\le\ell(N)$ and $K(i)\le2\log(i+2)+O(1)$, this gives
$\ell(N)\le 6\log(\ell(N)+2)+O(1)$. But $N\ge n_i$ and $\ell$ is unbounded, so $\ell(N)\to\infty$ as $i\to\infty$;
contradiction.

(b) The hypothesis and $K(x\rest n\mid (y\rest n)^*)\le n+O(1)$ give $K(y\rest n)\ge n-h(n)-O(1)$, so
$K(u\rest k)\ge k-h(k)-O(1)$ for all $k$. Apply Lemma~\ref{lem:slack} to $2h$ and argue as in (a):
since $g(N)\le N$ and $h$ is nondecreasing,
$K(u\rest g(N))\ge g(N)-h(N)-O(1)$, and we get $g(N)-2h(N)\le 3K(i)+O(1)$, the same contradiction.
\end{proof}

So if $x$ is compatible with a random $y$, then every computable $h$ with
$K(x\rest n,y\rest n)\ge 2n-h(n)-O(1)$ satisfies $\sum_n2^{-h(n)}<\infty$: the joint deficiency of a compatible pair
is \emph{summable}. Compare the ample excess lemma of \cite{MY}, by which $\sum_n2^{n-K(X\rest n)}<\infty$ for random $X$.

\section{Overwritten random reals}\label{sec:over}

For $F\subseteq\w$ and a real $R$, let $R_F$ be $R$ with the positions in $F$ overwritten by zeros:
$R_F(k)=0$ if $k\in F$ and $R_F(k)=R(k)$ otherwise. Let $f(n)=|F\cap[0,n)|$ be the counting function of $F$.
(If $F\subseteq G$ are computable then $R_G\leK R_F$, since $R_G\rest n$ is computable from $R_F\rest n$.
Nothing below changes if the positions in $F$ are overwritten by any other computable pattern.)

\begin{lemma}\label{lem:over}
Let $F$ be computable and let $R,Y$ be reals.
\begin{enumerate}[label=\textup{(\alph*)}]
\item $K(R_F\rest n)\le K(R\rest n)+O(1)$; so $R_F\leK R$.
\item $K(R\rest n,\,Y\rest n)\le K(R_F\rest n,\,Y\rest n)+f(n)+O(1)$.
\item $K(R_F\rest n)\le n-f(n)+K(n)+O(1)$.
\item If $F$ is infinite, $R_F$ is not random \textup{(}not even Kurtz random\textup{)}. Hence, if $R$ is random,
$R_F<_K R$.
\end{enumerate}
\end{lemma}
\begin{proof}
(a) $R_F\rest n$ is computable from $R\rest n$.
(b) From the pair one reads off $n$, hence $F\cap[0,n)$ and $f(n)$; the missing $f(n)$ bits of $R$ are then appended.
(c) Describe $n$, then the $n-f(n)$ bits of $R_F\rest n$ at positions outside $F$.
(d) $\{Z:Z(k)=0\text{ for all }k\in F\}$ is a null $\Pi^0_1$ class. If $R\leK R_F$, then $R_F$ would be random by \eqref{eq:up}.
\end{proof}

\begin{theorem}[Dichotomy]\label{thm:dich}
Let $F$ be a computable set with counting function $f$, and let $R\oplus Y$ be random.
\begin{enumerate}[label=\textup{(\roman*)}]
\item If $\sum_n2^{-f(n)}<\infty$, then $K(R_F\rest n)\le n+O(1)$, and so $R_F\leK W$ for \emph{every} random $W$.
In particular $R_F\leK Y$.
\item If $\sum_n2^{-f(n)}=\infty$, then $R_F$ and $Y$ have no common upper bound in the $K$-degrees.
\end{enumerate}
Thus $R_F$ and $Y$ are compatible iff $R_F\leK Y$ iff $\sum_n 2^{-f(n)}<\infty$.
\end{theorem}
\begin{proof}
(i) As $f$ is computable and $\sum_n2^{-f(n)}<\infty$, the Kraft--Chaitin theorem gives $K(n)\le f(n)+O(1)$.
By Lemma~\ref{lem:over}(c), $K(R_F\rest n)\le n+O(1)\le K(W\rest n)+O(1)$ for random $W$.
(This part needs neither the randomness of $R$ nor $Y$.)

(ii) By Lemma~\ref{lem:over}(b), \eqref{eq:join} and the randomness of $R\oplus Y$,
$K(R_F\rest n,Y\rest n)\ge 2n-f(n)-O(1)$ for all $n$. As $Y$ is random, Corollary~\ref{cor:crit}(a) with $h=f$ applies.
\end{proof}

\begin{corollary}\label{cor:J}
\begin{enumerate}[label=\textup{(\alph*)}]
\item For $F=\emptyset$, Theorem~\ref{thm:dich}\textup{(ii)} is the Miller--Yu theorem.
\item Let $F=\{2^k:k\in\w\}$ and let $R\oplus Y$ be random. Then $f(n)\le\log n+1$ for $n\ge1$, so $\sum_n2^{-f(n)}=\infty$.
Hence $R_F$ is not random, $R_F<_KR$, and $R_F\incomp Y$.
So $J$ holds of a non-random real: \emph{the formula $J$ does not define the random degrees}.
\item If $F$ is computable and $|F\cap[0,n)|\ge2\log n$ for all large $n$, then $R_F$ is $K$-below every random
real, for every real $R$.
\end{enumerate}
\end{corollary}

\begin{remark}[Two non-random reals with no common upper bound]
Let $F=\{2^{2^k}:k\in\w\}$ and let $R\oplus Y$ be random. Then
$K(R_F\rest n,Y_F\rest n)\ge 2n-2f(n)-O(1)$ with $2f(n)\le 2\log\log n+2$, and
$\sum_n2^{-4f(n)}=\infty$. By Corollary~\ref{cor:crit}(b) with $h=2f$, $R_F\incomp Y_F$, and neither real is random.
Section~\ref{sec:source} shows that $R_F\incomp Y_G$ whenever both counting functions are at most $\alpha\log n$
with $\alpha<1$.
\end{remark}

\section{Reals below every random real}\label{sec:other}

\begin{proposition}\label{prop:below}
If $K(x\rest n)\le n+O(1)$ for all $n$, then $x$ is $K$-below every random real, and therefore compatible with every
real that has a random upper bound. Consequently:
\begin{enumerate}[label=\textup{(\alph*)}]
\item if Question~9.8 has a positive answer, every such $x$ is compatible with every real;
\item a single pair $x\incomp y$ with $K(x\rest n)\le n+O(1)$ would refute Question~9.8, with $y$ as the counterexample.
\end{enumerate}
\end{proposition}
\begin{proof}
$K(x\rest n)\le n+O(1)\le K(w\rest n)+O(1)$ for random $w$. If $y\leK w$ with $w$ random, $w$ bounds $x$ and $y$.
\end{proof}

\begin{proposition}[using Zhao \cite{Zhao}]\label{prop:dim}
If $\limsup_nK(x\rest n)/n<1$, then $x$ is compatible with every real; that is, $J(x)$ fails.
\end{proposition}
\begin{proof}
By \cite[Theorem~1.1]{Zhao}, every real $y$ has an upper bound $B$ with $K(B\rest n)\ge n-o(n)$. Then
$K(x\rest n)\le K(B\rest n)$ for all large $n$, so $x\leK B$.
\end{proof}

Proposition~\ref{prop:dim} depends on the recent preprint \cite{Zhao}.

\medskip
\noindent\textbf{The picture for overwritten reals.}
Let $R\oplus Y$ be random and $F$ computable and infinite with counting function $f$. Then $R_F$ is non-random,
$R_F<_KR$, and:
\begin{center}
\begin{tabular}{ll}
$\sum_n2^{-f(n)}=\infty$ (e.g.\ $F=\{2^k\}$) & $R_F\incomp Y$; so $J(R_F)$ holds.\\[2pt]
$\sum_n2^{-f(n)}<\infty$ (e.g.\ $f(n)\ge2\log n$) & $R_F$ is below every random real;\\
 & $J(R_F)$ fails if Question~9.8 is true.\\[2pt]
$\liminf_n f(n)/n>0$ & $J(R_F)$ fails outright (Proposition~\ref{prop:dim}).
\end{tabular}
\end{center}
For random $R$ and computable $F$ we therefore have, unconditionally,
\[
 \sum_n2^{-f(n)}=\infty\iff\sup_n\bigl(K(R_F\rest n)-n\bigr)=\infty ,
\]
and both imply $J(R_F)$. The converse implication, from $J(R_F)$ to divergence, holds if Question~9.8 is true.

\section{What a random upper bound must do}\label{sec:need}

Applying the cover lemma to a random upper bound gives a necessary condition for Question~9.8 at a given real.

\begin{proposition}\label{prop:nec}
Suppose $x\leK w$ with $w$ random, and let $g$ be a gauge with blocks $I_i$. Then there is $c$ such that every block
$I_i$ contains a length $N$ with $g(N)<\infty$ and
\[
 K(x\rest N)\le N+K(i)+c\qquad\text{and}\qquad K(N)\ge g(N)-K(i)-c .
\]
\end{proposition}
\begin{proof}
Apply Lemma~\ref{lem:cover} to $u=w$. Part (a) and $x\leK w$ give the first inequality. Part (b) and
$K(w\rest g(N))\ge g(N)-O(1)$ give the second.
\end{proof}

Every real $x$ has, in every block, a length $N$ with $K(x\rest N)\le N+K(i)+O(1)$, namely the one that $x$ itself
announces. The content of Proposition~\ref{prop:nec} is that a random upper bound forces such a low point at a
length $N$ that is \emph{complex} relative to the gauge, $K(N)\ge g(N)-K(i)-O(1)$, because a random real can only
announce complex lengths. Hence:

\begin{corollary}\label{cor:route}
Suppose that for some gauge $g$ and every $c$ there is a block $I_i$ such that $K(x\rest N)>N+K(i)+c$ for every
$N\in I_i$ with $K(N)\ge g(N)-K(i)-c$. Then $x$ has no random $K$-upper bound.
\end{corollary}

We do not know whether any real satisfies the hypothesis of Corollary~\ref{cor:route}; such a real would refute
Question~9.8. It is a precise form of ``the low points of $x$ cannot be announced by a random real''.
The same argument gives a criterion for incompatibility with a fixed random $y$: it suffices that, for every $c$,
some block $I_i$ has $K(x\rest N)>N+K(i)+c$ at every $N\in I_i$ with $K(N)\ge g(N)-K(i)-c$ and
$K(y\rest N)\le N+K(i)+c$.

\begin{question}\label{q1}
Does $\sup_n\bigl(K(x\rest n)-n\bigr)=\infty$ imply $J(x)$, i.e.\ that $x$ is incompatible with some real?
\end{question}
The method above proves $x\incomp y$ only when $x$ is complex at essentially all lengths where $y$ can be low.
A real whose excess $K(x\rest n)-n$ is unbounded but positive only on a sparse set of lengths is not covered.

\begin{question}\label{q2}
Is every real with $K(x\rest n)\le n+O(1)$ compatible with every real? \textup{(}Yes if Question~9.8 is true; yes
when $\limsup_nK(x\rest n)/n<1$.\textup{)} Does a positive answer imply Question~9.8?
\end{question}
If Questions~\ref{q1} and~\ref{q2} both have positive answers, then the class
$\{x:K(x\rest n)\le n+O(1)\}$ is first-order definable in the $K$-degrees, as the set of degrees compatible with
every degree.

\begin{question}\label{q3}
Is a real random iff it is a $K$-upper bound of every real $z$ with $K(z\rest n)\le n+O(1)$?
\end{question}
With Questions~\ref{q1} and~\ref{q2} this would make randomness definable: $x$ is random iff $x$ bounds every degree
that is compatible with all degrees. The delicate case is a non-random $x$ of very slowly growing deficiency.

\begin{question}\label{q4}
For computable $F,G$ and random $R$, when is $R_F\leK R_G$, and when are $R_F$ and $R_G$ compatible?
\end{question}
Inclusion $G\subseteq F$ is sufficient for $R_F\leK R_G$. Theorem~\ref{thm:dich} treats $R_F$ against an independent
random real; Question~\ref{q4} asks for the corresponding analysis inside the lower cone of $R$, which is where a
coding of structures into the $K$-degrees would have to live. Proposition~\ref{prop:order} answers the first part
for sparse sets.

\section{Stage masks: which finite configurations occur}\label{sec:masks}

Sections \ref{sec:masks}--\ref{sec:profiles} concern the first-order theory of the $K$-degrees
$\mathcal D_K=(\cs/{\equiv_K},\le)$. This section determines the finite substructures of $(\mathcal D_K,\le,\incomp)$.

Let $N_0=1$ and $N_{t+1}=2^{8N_t}$. For $t\ge1$, \emph{stage} $t$ is the interval $P_t=[N_{t-1},N_t)$.
Fix reals $\rho^1,\rho^2$ with $\rho^1\oplus\rho^2$ random.
A \emph{pattern} is a computable function $M:\{1,2,\dots\}\to\{0,1,2\}$, and
\[
 x_M(k)=\begin{cases}\rho^{M(t)}(k)&\text{if }k\in P_t\text{ and }M(t)\ne0,\\ 0&\text{otherwise.}\end{cases}
\]
So on stage $t$ the real $x_M$ copies the source $\rho^{M(t)}$, or is inactive if $M(t)=0$. For patterns $M,M'$ write
\begin{align*}
M\le^*M' &\iff M(t)\in\{0,M'(t)\}\text{ for all but finitely many }t;\\
M\incomp^*M' &\iff M(t)\text{ and }M'(t)\text{ are distinct and nonzero for infinitely many }t.
\end{align*}

\begin{lemma}\label{lem:seg}
If $Z$ is random, then $K\bigl(Z(s)Z(s+1)\cdots Z(m-1)\bigr)\ge m-s-K(s)-O(1)$ for all $s\le m$.
\end{lemma}
\begin{proof}
Concatenating programs, $K(Z\rest m)\le K(Z\rest s)+K(Z(s)\cdots Z(m-1))+O(1)$. Now use $K(Z\rest s)\le s+K(s)+O(1)$ and
$K(Z\rest m)\ge m-O(1)$.
\end{proof}

\begin{theorem}\label{thm:masks}
For all patterns $M,M'$:
\begin{enumerate}[label=\textup{(\alph*)}]
\item $x_M\leK x_{M'}\iff M\le^*M'$;
\item $x_M\incomp x_{M'}\iff M\incomp^*M'$.
\end{enumerate}
\end{theorem}
\begin{proof}
\emph{Step 1: if $M\le^*M'$ then $x_M\leK x_{M'}$.} With finitely many bits hard-wired, $x_M\rest n$ is computable
from $x_{M'}\rest n$: on a stage with $M(t)=M'(t)$ copy, on a stage with $M(t)=0$ write zeros.

\emph{Step 2: if $M(t)\ne0=M'(t)$ for infinitely many $t$, then $x_M\not\leK x_{M'}$.} For such $t$, $x_M\rest N_t$
computes the bits of a random real on $P_t$, so $K(x_M\rest N_t)\ge N_t-N_{t-1}-K(N_{t-1})-O(1)$ by Lemma~\ref{lem:seg}.
But $x_{M'}\rest N_t$ is $x_{M'}\rest N_{t-1}$ followed by zeros, so
$K(x_{M'}\rest N_t)\le N_{t-1}+K(N_{t-1})+K(t)+O(1)$. The difference is unbounded.

\emph{Step 3: if $M\incomp^*M'$ then $x_M\incomp x_{M'}$.} Let $T$ be the infinite computable set of $t$ at which
$M(t),M'(t)$ are distinct and nonzero, and suppose $u$ is a common upper bound. Define
\[
 g(N)=8N_{t-1}\ \text{ if }N\in[8N_{t-1},N_t)\text{ for some }t\in T,\qquad g(N)=\infty\ \text{ otherwise.}
\]
Then $g(N)\le N$ where finite, and each stage $t\in T$ contributes
$\sum_{N\in P_t}2^{-g(N)}=1-8N_{t-1}/N_t\in[\tfrac12,1)$. So $g$ is a gauge, and every block of $g$ meets at least
two stages of $T$; hence if $N\in I_i$ lies in stage $t$ then $i\le t$.
By Theorem~\ref{thm:gauge}, for each $i$ there are $t\in T$ with $t\ge i$ and $N\in I_i\cap[8N_{t-1},N_t)$ with
\begin{equation}\label{eq:maskup}
 K(x_M\rest N,\;x_{M'}\rest N)\le 2N-K(u\rest 8N_{t-1})+3K(i)+O(1).
\end{equation}
On stage $t$ one of $x_M,x_{M'}$ copies $\rho^1$ and the other copies $\rho^2$. So the pair
$(x_M\rest N,x_{M'}\rest N)$ computes the bits of $\rho^1\oplus\rho^2$ at the positions in $[2N_{t-1},2N)$, and by
Lemma~\ref{lem:seg} the left side of \eqref{eq:maskup} is at least $2N-2N_{t-1}-K(2N_{t-1})-O(1)$.
Likewise $x_M\rest 8N_{t-1}$ computes the bits of a random real on $[N_{t-1},8N_{t-1})$, so
$K(u\rest 8N_{t-1})\ge K(x_M\rest 8N_{t-1})-O(1)\ge 7N_{t-1}-K(N_{t-1})-O(1)$.
Substituting, $5N_{t-1}\le 2K(N_{t-1})+3K(i)+O(1)\le 4\log N_{t-1}+6\log(t+2)+O(1)$, which is false for
large $t$ since $N_{t-1}\ge t$. As $t\ge i$, this is a contradiction for large $i$.

\emph{Step 4: if not $M\incomp^*M'$ then $x_M,x_{M'}$ are compatible.} Let $M''(t)=M(t)$ if $M(t)\ne0$ and
$M''(t)=M'(t)$ otherwise. Then $M\le^*M''$ and $M'\le^*M''$, so $x_{M''}$ is a common upper bound by Step~1.

Now (b) is Steps 3 and 4. For (a), one direction is Step 1. If $M\not\le^*M'$ there are infinitely many $t$ with
$M(t)\ne0$ and $M(t)\ne M'(t)$; either infinitely many of them have $M'(t)=0$ and Step~2 applies, or $M\incomp^*M'$
and by Step~3 the two reals are not even compatible.
\end{proof}

Call $(P,\le,\incomp)$ a \emph{compatibility structure} if $\le$ is a partial order on $P$, $\incomp$ is a symmetric
relation on $P$, and for all $a,a',b\in P$:
\begin{enumerate}[label=\textup{(C\arabic*)}]
\item $a\le b$ implies not $a\incomp b$;
\item $a\incomp b$ and $a\le a'$ imply $a'\incomp b$.
\end{enumerate}
$(\mathcal D_K,\le,\incomp)$ is a compatibility structure (Proposition~\ref{prop:basic}(a)).

\begin{corollary}\label{cor:fin}
\begin{enumerate}[label=\textup{(\alph*)}]
\item For every random $\rho$, the map $A\mapsto$ ``$\rho$ masked to the stages in $A$'' embeds the Boolean algebra of
computable sets modulo finite sets, as a partial order, into the $K$-degrees below $\rho$.
\item Every finite compatibility structure is isomorphic to a substructure of $(\mathcal D_K,\le,\incomp)$.
In particular every finite partial order embeds in $\mathcal D_K$, and every finite graph is the compatibility graph
of a finite antichain of $K$-degrees.
\item The existential theory of $(\mathcal D_K,\le,\incomp)$ is decidable.
\end{enumerate}
\end{corollary}
\begin{proof}
(a) Use the patterns with values in $\{0,1\}$ and $\rho^1=\rho$; then $M\le^*M'$ is inclusion modulo finite sets.
Only Steps 1 and 2 are used, so one random source suffices; and $x_M\leK\rho$ because $x_M\rest n$ is computable
from $\rho\rest n$.

(b) Let $(P,\le,\incomp)$ be finite. Partition $\{1,2,\dots\}$ into infinite computable sets $T_d$ ($d\in P$) and
$T_{e,e'}$, one for each unordered pair with $e\incomp e'$, listed in a fixed order $(e,e')$. For $a\in P$ let
\[
M_a(t)=\begin{cases}
1&\text{if }t\in T_d\text{ and }d\le a,\ \text{ or }t\in T_{e,e'}\text{ and }e\le a;\\
2&\text{if }t\in T_{e,e'}\text{ and }e'\le a;\\
0&\text{otherwise.}\end{cases}
\]
This is well defined: by (C2) and (C1) no $a$ is above both $e$ and $e'$ when $e\incomp e'$.
If $a\le b$ then $M_a(t)\in\{0,M_b(t)\}$ for all $t$. If $a\not\le b$ then $M_a=1$ and $M_b=0$ on $T_a$, so
$M_a\not\le^*M_b$. If $a\incomp b$ then $M_a,M_b$ take the values $1,2$ on $T_{a,b}$ (or on $T_{b,a}$), so $M_a\incomp^*M_b$.
If $M_a(t),M_b(t)$ are distinct and nonzero then $t\in T_{e,e'}$ with, say, $e\le a$ and $e'\le b$, and two
applications of (C2) give $a\incomp b$. By Theorem~\ref{thm:masks}, $a\mapsto x_{M_a}$ is the required embedding.
For graphs take $P$ an antichain and $\incomp$ the set of non-edges.

(c) By (b), an existential sentence holds in $(\mathcal D_K,\le,\incomp)$ iff its matrix is satisfiable in some
finite compatibility structure with at most as many elements as the sentence has variables.
\end{proof}

Since $\incomp$ is defined by a universal formula of the order language and compatibility by an existential one,
Corollary~\ref{cor:fin}(c) decides a fragment of the $\exists\forall$-theory of $(\mathcal D_K,\le)$, and in particular
its existential theory.

\begin{remark}[Every ideal is cut out by one degree]\label{rem:ideal}
For an arbitrary set $D$ of stages let $x_D$ be $\rho$ masked to $D$. If $A$ is computable then
$x_A\leK x_D$ iff $A\subseteq^*D$: the proofs of Steps 1 and 2 do not use that the larger mask is computable.
Let $I$ be any ideal of the Boolean algebra of computable sets modulo finite sets, $I=\{A_0,A_1,\dots\}$, and let
$C_0,C_1,\dots$ list the computable sets not in $I$. Choose $m_0<m_1<\cdots$ such that
$\bigl(C_j\setminus(A_0\cup\dots\cup A_k)\bigr)\cap[m_k,m_{k+1})\ne\emptyset$ for all $j\le k$, and put
$D=\bigcup_k(A_k\setminus[0,m_k))$. Then for computable $A$: $A\subseteq^*D$ iff $A\in I$. So
$\{A: x_A\leK x_D\}=I$.

Thus below every random degree there is a copy of the countable atomless Boolean algebra all of whose ideals are
traces of lower cones. By itself this does not give undecidability: by a theorem of Rabin \cite{Rabin} the theory
of countable Boolean algebras with quantification over ideals is decidable. But it is the exactness property on
which codings with parameters rest.
\end{remark}

\begin{lemma}[Long common blocks]\label{lem:blocks}
Let $x,y$ be reals, and suppose that for infinitely many $t$ we have $x(k)=\rho^1(k)$ and $y(k)=\rho^2(k)$ for all
$k$ with $N_{t-1}\le k<2^{3N_{t-1}+1}$. Then $x\incomp y$.
\end{lemma}
\begin{proof}
Suppose $x,y\leK u$, and let $c$ be as in Theorem~\ref{thm:basic}. Fix such a $t$ and put $m=3N_{t-1}$ and $N=u_m$,
so $N_{t-1}\le2^m\le N<2^{m+1}$. By Theorem~\ref{thm:basic}, $K(x\rest N,y\rest N)\le 2N-m+K(m)+2c$. On the other
hand $(\rho^1\oplus\rho^2)\rest 2N$ is computable from $(\rho^1\oplus\rho^2)\rest 2N_{t-1}$ and the pair
$(x\rest N,y\rest N)$, so
$2N-O(1)\le K((\rho^1\oplus\rho^2)\rest 2N)\le 2N_{t-1}+K(2N_{t-1})+K(x\rest N,y\rest N)+O(1)$. Hence
$N_{t-1}\le K(3N_{t-1})+K(2N_{t-1})+O(1)\le 4\log N_{t-1}+O(1)$, which is false for large $t$.
\end{proof}
Lemma~\ref{lem:blocks} gives Step~3 of Theorem~\ref{thm:masks} without gauges, and it does not require the two
reals to be given by computable patterns: only the first $2^{3N_{t-1}+1}$ positions of infinitely many common
stages matter. (Lemma~\ref{lem:short} shortens the blocks to $2^{N+O(\log N)}$.)

\section{Countable lower cones}\label{sec:sparse}

A random real has a large lower cone: every real of bounded excess is below it (Proposition~\ref{prop:below}).
Reals with long runs of zeros behave in the opposite way.

\begin{definition}
A length $\ell$ is \emph{announced} in a real $x$, with constant $e$, if there is $n\ge\ell$ such that
$x(k)=0$ for all $\ell\le k<n$, and $K(n\mid(x\rest\ell)^*)\le e$ and $K(x\rest\ell\mid n)\le e$.
\end{definition}
For instance, if $x(k)=0$ for all $\ell\le k<2^{\ell+1}$, then $\ell$ is announced, with a constant not depending on
$x$ or $\ell$: take $n$ to be the integer with binary expansion $1\,x(0)\cdots x(\ell-1)$.

\begin{proposition}[Announced prefixes are rigid]\label{prop:rigid}
Suppose $\ell$ is announced in $x$ with constant $e$, and $K(y\rest n)\le K(x\rest n)+c$ for all $n$. Then
\[
 K\bigl(y\rest\ell\mid(x\rest\ell)^*\bigr)\le c+2e+O(1).
\]
\end{proposition}
\begin{proof}
Let $n$ witness that $\ell$ is announced. Since $x\rest n$ is $x\rest\ell$ followed by zeros,
$K(x\rest n)\le K(x\rest\ell)+K(n\mid(x\rest\ell)^*)+O(1)\le K(x\rest\ell)+e+O(1)$, so
$K(y\rest n)\le K(x\rest\ell)+e+c+O(1)$. By \eqref{eq:soi},
\begin{align*}
K(y\rest n\mid(x\rest\ell)^*)&=K(x\rest\ell,\,y\rest n)-K(x\rest\ell)+O(1)\\
&\le K(y\rest n)+K(x\rest\ell\mid n)-K(x\rest\ell)+O(1)\ \le\ c+2e+O(1).
\end{align*}
Finally $y\rest\ell$ is computable from $y\rest n$ and $\ell$, and $\ell$ is known from $(x\rest\ell)^*$.
\end{proof}

\begin{corollary}\label{cor:countable}
If infinitely many lengths are announced in $x$ with one constant, then only countably many reals are $K$-below
$x$. This holds in particular
\begin{enumerate}[label=\textup{(\alph*)}]
\item if $x$ vanishes on $[\ell,2^{\ell+1})$ for infinitely many $\ell$; for example, if $x$ is any real masked to a
co-infinite set $A$ of stages \textup{(}in the sense of Section~\ref{sec:masks}, with $\rho$ and $A$ otherwise
arbitrary\textup{)};
\item for the sparse codes $S_V$ defined below.
\end{enumerate}
\end{corollary}
\begin{proof}
For each $d$ and each announced $\ell$, fewer than $2^{d+1}$ strings $\sigma$ of length $\ell$ satisfy
$K(\sigma\mid(x\rest\ell)^*)\le d$. So by Proposition~\ref{prop:rigid}, for each $c$ only finitely many reals $y$
satisfy $K(y\rest n)\le K(x\rest n)+c$ for all $n$. For (a) with a mask, take $\ell=N_t$ with $t+1\notin A$: the mask vanishes on
$[N_t,N_{t+1})\supseteq[\ell,2^{\ell+1})$. For (b) see the proof of Proposition~\ref{prop:sparse}.
\end{proof}

Corollary~\ref{cor:countable} is a special case of a known fact. A real $x$ is \emph{infinitely often $K$-trivial} if
$K(x\rest n)\le K(n)+O(1)$ for infinitely many $n$, and Barmpalias and Vlek \cite{BV} showed, by the same counting,
that such a real has only countably many predecessors; see also \cite{BL2}. If $\ell$ is announced in $x$ with
witness $n$, then $K(x\rest n)\le K(n)+O(1)$, so the reals of the corollary are infinitely often $K$-trivial. What
Proposition~\ref{prop:rigid} adds is the form of the predecessors: at an announced length they are computable,
with bounded advice, from a shortest program for the prefix of $x$.

For a real $V$ let $S_V$ be the real with $S_V(4^t)=V(t)$ for all $t$ and $S_V(k)=0$ when $k$ is not a power of $4$.
For these sparse codes the lower cone can be described exactly.

\begin{proposition}\label{prop:sparse}
For all reals $y$ and $V$:
\[
 y\leK S_V\iff \sup_t K\bigl(y\rest 4^t\mid (V\rest t)^*\bigr)<\infty .
\]
\end{proposition}
\begin{proof}
For $4^{t-1}<n\le4^t$ the string $S_V\rest n$ and the pair $(n,V\rest t)$ are computable from each other, so
$K(S_V\rest n)=K(n,V\rest t)\pm O(1)$.

($\Leftarrow$) Suppose $K(y\rest 4^t\mid(V\rest t)^*)\le d$ for all $t$, and let $4^{t-1}<n\le4^t$. Since $y\rest n$
is a prefix of $y\rest4^t$, $K(y\rest n)\le K\bigl(n,(V\rest t)^*\bigr)+d+O(1)$. By \eqref{eq:soi},
$K(n,(V\rest t)^*)=K(V\rest t)+K(n\mid(V\rest t)^*)+O(1)=K(n,V\rest t)+O(1)$. So $K(y\rest n)\le K(S_V\rest n)+d+O(1)$.

($\Rightarrow$) Suppose $K(y\rest n)\le K(S_V\rest n)+c$ for all $n$. Fix $t$, let $\tau=V\rest(t+1)$ and $\ell=4^t+1$, and let $n=4^t+1+\mathrm{num}(\tau)$, where
$\mathrm{num}(\tau)<2^{t+1}$ is the number with binary expansion $\tau$. Then $\ell\le n\le4^{t+1}$, $S_V$ vanishes on
$[\ell,n)$, and the three objects $n$, $\tau$ and $S_V\rest\ell$ are computable from one another. So $\ell$ is
announced in $S_V$ with an absolute constant, and Proposition~\ref{prop:rigid} gives
$K(y\rest4^t\mid(S_V\rest\ell)^*)\le c+O(1)$. Finally $(S_V\rest\ell)^*$ is
computable from $(V\rest t)^*$ and $O(1)$ further bits, because $K(\tau)=K(V\rest t)\pm O(1)$.
\end{proof}

\begin{corollary}\label{cor:sparse}
$S_W\leK S_V$ iff $\sup_tK(W\rest t\mid(V\rest t)^*)<\infty$. In particular $W\le_{rK}V$ implies $S_W\leK S_V$.
\end{corollary}
\begin{proof}
$S_W\rest4^t$ and $W\rest t$ are computable from each other given $t$.
\end{proof}

So on sparse codes $\leK$ is a strong reducibility between the coded reals: a slight weakening of
$rK$-reducibility ($K(W\rest t\mid V\rest t)\le O(1)$), in which the condition is a shortest program for $V\rest t$.
The lower cone of $S_V$ consists exactly of the reals whose prefixes of length $4^t$ are uniformly simple given a
shortest program for
$V\rest t$. This is the kind of setting in which initial-segment constructions are usually possible.
(The existence of a minimal $K$-degree is nevertheless open; see Section~\ref{sec:fo}.)

\section{Joins and meets}\label{sec:timing}

The coding discussed in Section~\ref{sec:fo} needs a Boolean algebra whose operations are visible to the order of
$\mathcal D_K$. Least upper bounds need not exist in $\mathcal D_K$ \cite{MY}, so this is a real constraint.
This section gives a family of reals that meets it, and shows that the stage masks of Section~\ref{sec:masks} do not.

\begin{lemma}\label{lem:maxmin}
Let $x,y,u,w$ be reals.
\begin{enumerate}[label=\textup{(\alph*)}]
\item If $K(u\rest n)=\max\{K(x\rest n),K(y\rest n)\}\pm O(1)$ for all $n$, then the degree of $u$ is the least upper
bound of the degrees of $x$ and $y$ in $\mathcal D_K$.
\item If $K(w\rest n)=\min\{K(x\rest n),K(y\rest n)\}\pm O(1)$ for all $n$, then the degree of $w$ is their greatest
lower bound in $\mathcal D_K$.
\end{enumerate}
\end{lemma}
\begin{proof}
(a) Clearly $x,y\leK u$. If $x,y\leK v$ with constant $c$, then for all $n$
\[
 K(u\rest n)\le\max\{K(x\rest n),K(y\rest n)\}+O(1)\le K(v\rest n)+c+O(1).
\]
(b) is dual.
\end{proof}

So joins and meets exist whenever the pointwise maximum and minimum of the two complexity functions are again
complexity functions of reals. Lemma~\ref{lem:maxmin} itself does not say when this happens; for the timing masks
below, the reals that realize the maximum and the minimum are constructed explicitly.

\smallskip
\noindent\textbf{Timing masks.}
Keep $N_t$ and the stages $P_t=[N_{t-1},N_t)$ of Section~\ref{sec:masks}. Put $\ell_t=N_{t-1}$ and
$c_t=\ell_1+\dots+\ell_t$, with $c_0=0$; note that $2\ell_t+1\le N_t-N_{t-1}$. Fix a random real $\rho$ and let
$\Gamma_t=\rho(c_{t-1})\cdots\rho(c_t-1)$ be its $t$-th block, of length $\ell_t$. For a set $S\subseteq\{1,2,\dots\}$ let
$z_S$ be the real with $z_S(0)=0$ which on stage $t$ is defined as follows:
\begin{itemize}[leftmargin=2.2em]
\item $z_S(N_{t-1})=S(t)$ (a flag);
\item if $t\in S$, the block $\Gamma_t$ is written on the positions $N_{t-1}+1,\dots,N_{t-1}+\ell_t$ (\emph{early});
\item if $t\notin S$, the block $\Gamma_t$ is written on the last $\ell_t$ positions of $P_t$ (\emph{late});
\item all other bits of $z_S$ on $P_t$ are $0$.
\end{itemize}
All the reals $z_S$ carry the same information, namely $\rho$. The set $S$ only decides \emph{when} in each stage the
next block appears. This is the difference from Section~\ref{sec:masks}, where a mask decides \emph{which} blocks
appear. For $N_{t-1}<n\le N_t$ let
\[
 m_S(n)=\begin{cases} c_{t-1}+\min\{n-N_{t-1}-1,\ \ell_t\}&\text{if }t\in S,\\
 c_{t-1}+\max\{0,\ n-(N_t-\ell_t)\}&\text{if }t\notin S,\end{cases}
\]
the number of bits of $\rho$ that $z_S$ has shown below $n$.

\begin{lemma}\label{lem:level}
\begin{enumerate}[label=\textup{(\alph*)}]
\item For every set $D$, the string $\rho\rest m_D(n)$ is computable from $z_D\rest n$, uniformly in $n$ and $D$.
\item If $S$ is computable, then $z_S\rest n$ and the pair $(n,\rho\rest m_S(n))$ are computable from each other,
uniformly in $n$. Hence $K(z_S\rest n)=K(n,\rho\rest m_S(n))\pm O(1)$.
\item $m_{S\cup S'}=\max\{m_S,m_{S'}\}$ and $m_{S\cap S'}=\min\{m_S,m_{S'}\}$.
\end{enumerate}
\end{lemma}
\begin{proof}
(a) The flags below $n$ say where the blocks sit. (b) In addition, $z_S\rest n$ is determined by $n$, $S$ and the
bits of $\rho$ it shows. (c) On each stage the early schedule has shown at least as much as the late one:
$\min\{n-N_{t-1}-1,\ell_t\}\ge\max\{0,n-(N_t-\ell_t)\}$ because $N_t-N_{t-1}\ge2\ell_t+1$.
\end{proof}

\begin{theorem}\label{thm:timing}
Let $S,S'$ be computable sets.
\begin{enumerate}[label=\textup{(\alph*)}]
\item $z_S\leK z_{S'}$ iff $S\subseteq^*S'$.
\item The degree of $z_{S\cup S'}$ is the least upper bound, and the degree of $z_{S\cap S'}$ the greatest lower
bound, of the degrees of $z_S$ and $z_{S'}$ in $\mathcal D_K$.
\end{enumerate}
\end{theorem}
\begin{proof}
(b) Fix $n$ and suppose $m_S(n)\le m_{S'}(n)$; the other case is symmetric. By Lemma~\ref{lem:level}(c),
$m_{S\cup S'}(n)=m_{S'}(n)$ and $m_{S\cap S'}(n)=m_S(n)$. So by Lemma~\ref{lem:level}(b),
\[
 K(z_{S\cup S'}\rest n)=K(z_{S'}\rest n)\pm O(1)\quad\text{and}\quad K(z_{S\cap S'}\rest n)=K(z_S\rest n)\pm O(1).
\]
Moreover the pair $(n,\rho\rest m_S(n))$ is computable from the pair $(n,\rho\rest m_{S'}(n))$, so
$K(z_S\rest n)\le K(z_{S'}\rest n)+O(1)$.
Hence, with constants that do not depend on $n$,
\begin{align*}
K(z_{S\cup S'}\rest n)&=\max\{K(z_S\rest n),K(z_{S'}\rest n)\}\pm O(1),\\
K(z_{S\cap S'}\rest n)&=\min\{K(z_S\rest n),K(z_{S'}\rest n)\}\pm O(1).
\end{align*}
Now apply Lemma~\ref{lem:maxmin}.

(a) If $S\subseteq^*S'$ then $m_S(n)\le m_{S'}(n)$ for all sufficiently large $n$, and for these $n$ we have just
seen that $K(z_S\rest n)\le K(z_{S'}\rest n)+O(1)$. Conversely, suppose $S\setminus S'$ is infinite. For
$t\in S\setminus S'$ let $n_t=N_{t-1}+1+\ell_t$. Then $m_S(n_t)=c_t$, and $m_{S'}(n_t)=c_{t-1}$ because
$n_t\le N_t-\ell_t$. By Lemma~\ref{lem:level}(b), $K(z_S\rest n_t)\ge K(\rho\rest c_t)-O(1)\ge c_t-O(1)$ and
$K(z_{S'}\rest n_t)\le K(n_t,\rho\rest c_{t-1})+O(1)\le c_{t-1}+K(t)+O(1)$. The difference is at least
$\ell_t-K(t)-O(1)=N_{t-1}-K(t)-O(1)$, which is unbounded.
\end{proof}

For reals $a\leK b$ write $[a,b]$ for the set of $K$-degrees between them, and let $\mathrm{Comp}[a,b]$ be the set of
degrees $y\in[a,b]$ for which there is $y'\in[a,b]$ such that $a$ is the greatest lower bound and $b$ the least
upper bound of $y$ and $y'$ in $\mathcal D_K$. This set is first-order definable in $\mathcal D_K$ from the parameters
$a$ and $b$.

\begin{corollary}\label{cor:comp}
Let $a=z_\emptyset$ and $b=z_{\w}$, where $\w$ stands for the set of all stages. The map $S\mapsto z_S$ embeds the
Boolean algebra of computable sets modulo finite sets into the interval $[a,b]$ of $\mathcal D_K$ in such a way that
unions and intersections become least upper bounds and greatest lower bounds in $\mathcal D_K$. In particular
$z_S\in\mathrm{Comp}[a,b]$ for every computable $S$, with $z_{\w\setminus S}$ as a complement.
\end{corollary}

\begin{remark}\label{rem:timing}
(1) \emph{$K$-trivial sets.} Lemma~\ref{lem:level}(b), Theorem~\ref{thm:timing} and Corollary~\ref{cor:comp} hold for
$K$-trivial sets $S,S'$ as well. Such sets are low for $K$, that is $K(\sigma)=K^S(\sigma)\pm O(1)$
\cite[Ch.~11]{DH}, and relative to $S$ the two objects in Lemma~\ref{lem:level}(b) are computable from each other,
so $K(z_S\rest n)=K^S(z_S\rest n)\pm O(1)=K^S(n,\rho\rest m_S(n))\pm O(1)=K(n,\rho\rest m_S(n))\pm O(1)$. The
$K$-trivial sets are closed under union, intersection and complement.

(2) \emph{Names for ideals.} If $S$ is computable and $D$ is an arbitrary set, then $z_S\leK z_D$ iff
$S\subseteq^*D$. Indeed, if $S\subseteq^*D$ then $m_S(n)\le m_D(n)$ for large $n$, and by
Lemma~\ref{lem:level}(a) the pair $(n,\rho\rest m_S(n))$ is computable from $z_D\rest n$. If $t\in S\setminus D$
then $z_D\rest n_t$ is determined by $t$, the bits $D(1),\dots,D(t)$ and $\rho\rest c_{t-1}$, so
$K(z_D\rest n_t)\le c_{t-1}+t+K(t)+O(1)$, while $K(z_S\rest n_t)\ge c_t-O(1)$. Hence, exactly as in
Remark~\ref{rem:ideal}, every ideal of the algebra of computable sets modulo finite sets has the form
$\{S:z_S\leK z_D\}$ for some set $D$.
\end{remark}

The stage masks of Section~\ref{sec:masks} also embed the computable sets modulo finite sets
(Corollary~\ref{cor:fin}(a)), but there the union of two sets is not the join of the two degrees, not even inside
the lower cone of the mask of the union.

\begin{proposition}\label{prop:nojoin}
Let $\rho$ be random, let $A$ be the set of even numbers and $E$ the set of multiples of $4$, and let
$x_A,x_E,x_{A\setminus E}$ be $\rho$ masked to the stages in these sets, as in Section~\ref{sec:masks}. Then there
is a real $u$ with
\[
 x_E\leK u,\qquad x_{A\setminus E}\leK u,\qquad u\leK x_A,\qquad x_A\not\leK u .
\]
So $x_A$ is not a minimal upper bound of $x_E$ and $x_{A\setminus E}$.
\end{proposition}
\begin{proof}
For $t\equiv2\pmod 4$, $t\ge6$, let $D_t=N_{t-3}$. Let $u$ be $x_A$, except that for each such $t$ the last $D_t$
bits of stage $t$ are postponed by $D_t$ positions: $u(k)=0$ for $N_t-D_t\le k<N_t$, and
$u(N_t+i)=\rho(N_t-D_t+i)$ for $i<D_t$. Stage $t+1$ is odd, so $x_A$ vanishes on $[N_t,N_t+D_t)$.

\emph{$u\leK x_A$ and $x_E\leK u$.} $u\rest n$ is computable from $x_A\rest n$, and $x_E\rest n$ is computable
from $u\rest n$ because the stages in $E$ are untouched.

\emph{$x_A\not\leK u$.} Let $t\equiv2\pmod4$ and let $\alpha_t$ be the number of positions below $N_t-D_t$ that lie
in stages of $A$. The string $u\rest N_t$ is determined by $t$ and the bits of $\rho$ at these positions, so
$K(u\rest N_t)\le\alpha_t+K(t)+O(1)$. On the other hand $\rho\rest N_t$ is determined by $x_A\rest N_t$ and the bits
of $\rho$ at the $N_t-\alpha_t-D_t$ positions below $N_t$ outside the stages of $A$, so
$N_t-O(1)\le K(\rho\rest N_t)\le K(x_A\rest N_t)+N_t-\alpha_t-D_t+O(1)$. Hence
$K(x_A\rest N_t)-K(u\rest N_t)\ge D_t-K(t)-O(1)$, which is unbounded.

\emph{$x_{A\setminus E}\leK u$.} If $n$ lies in no interval $(N_t-D_t,N_t+D_t)$ then $x_{A\setminus E}\rest n$ is
computable from $u\rest n$: undo the postponements, which are complete below $n$, and erase the stages in $E$.
Now fix $t\equiv2\pmod4$, $t\ge6$, and $n$ with $N_t-D_t<n<N_t+D_t$. Put $m=N_t-D_t$ and
$Y=x_{A\setminus E}\rest m$, and let $\Gamma$ be the string of bits of $\rho$ on stage $t-2$, which is a stage of
$E$; so $|\Gamma|=N_{t-2}-N_{t-3}$. Both $Y$ and $\Gamma$ are computable from $u\rest m$, hence
\[
 K(u\rest n)\ge K(Y,\Gamma)-O(1).
\]
The string $x_{A\setminus E}\rest n$ is $Y$ followed by at most $D_t$ bits of $\rho$ and then by fewer than $D_t$
zeros, so
\[
 K(x_{A\setminus E}\rest n)\le K(Y)+D_t+4\log(D_t+1)+O(1).
\]
We claim that $K(Y,\Gamma)\ge K(Y)+|\Gamma|-K(t)-O(1)$. The positions below $m$ fall into three sets: stage $t-2$;
the positions in stages $\equiv2\pmod4$, say $\beta$ of them; and a set $R$ of remaining positions. So
$m=|\Gamma|+\beta+|R|$. The string $\rho\rest m$ is determined by $Y$, $\Gamma$ and the bits of $\rho$ on $R$, so
$m-O(1)\le K(\rho\rest m)\le K(Y,\Gamma)+|R|+O(1)$; and $Y$ is determined by $t$ and $\beta$ bits of $\rho$, so
$K(Y)\le\beta+K(t)+O(1)$. Subtracting gives the claim. Since $|\Gamma|=2^{8N_{t-3}}-N_{t-3}$, $D_t=N_{t-3}$ and
$t\le N_{t-3}+3$, we have $|\Gamma|-K(t)\ge D_t+4\log(D_t+1)$, and therefore
$K(x_{A\setminus E}\rest n)\le K(Y,\Gamma)+O(1)\le K(u\rest n)+O(1)$.
\end{proof}

The reason is that near the end of a stage of $A\setminus E$ the real $u$ still carries the earlier block of
$E$, which $x_{A\setminus E}$ lacks, and this surplus pays for the postponed bits. Nothing of the kind is possible
for timing masks, because there one of the two reals already has, at every length, all the information of the
other.

\begin{corollary}\label{cor:branch}
Let $S,S'$ be computable sets such that $S\setminus S'$ and $S'\setminus S$ are infinite. Then the degrees of $z_S$
and $z_{S'}$ are incomparable, their greatest lower bound in $\mathcal D_K$ is the degree of $z_{S\cap S'}$, and their
least upper bound is the degree of $z_{S\cup S'}$. In particular the degree of $z_{S\cap S'}$ is a non-zero
\emph{branching} degree: it is the greatest lower bound of two degrees strictly above it.
The same statements hold in the $C$-degrees, where $X\le_C Y$ means $C(X\rest n)\le C(Y\rest n)+O(1)$ for plain
complexity $C$.
\end{corollary}
\begin{proof}
For $\mathcal D_K$ this is Theorem~\ref{thm:timing}. The degree of $z_{S\cap S'}$ is not the least degree, because
$K(z_{S\cap S'}\rest N_t)\ge c_t-O(1)$ while $K(N_t)\le K(t)+O(1)$.
For plain complexity, Lemma~\ref{lem:maxmin}, Lemma~\ref{lem:level} and Theorem~\ref{thm:timing}(b) hold with $C$ in
place of $K$, with the same proofs: they use only that a computable map increases complexity by at most a
constant. In the proof of Theorem~\ref{thm:timing}(a) use $C(\rho\rest c_t)\ge c_t-2\log c_t-O(1)$, which follows
from $K(\sigma)\le C(\sigma)+2\log|\sigma|+O(1)$, and $C(n_t,\rho\rest c_{t-1})\le c_{t-1}+2\log t+O(1)$. Since
$c_t\le2N_{t-1}$, the difference is at least $N_{t-1}-2\log N_{t-1}-2\log t-O(1)$, which is unbounded. Finally
$z_{S\cap S'}$ computes $\rho$, so it is not computable, and the least $C$-degree consists of the computable sets.
\end{proof}

Miller and Nies asked \cite[Question~9.1]{MN} whether meets always exist in the vL-, $K$- and $C$-degrees, and
whether there is a non-zero branching degree. Merkle and Stephan \cite{MS} answered both parts for the $C$- and
$K$-degrees: some pairs have no meet, and some non-zero degrees are the meet of two other degrees.
Corollary~\ref{cor:branch} gives explicit examples in which the meet and the join both exist, and
Theorem~\ref{thm:timing} realizes them coherently on a whole Boolean algebra. Merkle and Stephan also embed the
lattice of subsets of $\omega$ \cite[Theorem~20]{MS} and every countable distributive lattice
\cite[Theorem~23]{MS} into these degree structures, with meets and joins preserved in the whole structure. Their proof of
\cite[Theorem~20]{MS} already realizes joins and meets as pointwise maxima and minima of complexity functions. So
Theorem~\ref{thm:timing} is not new as an embedding result, and the pointwise method is theirs. What is specific to
Theorem~\ref{thm:timing} is the construction: all reals carry one common random source, revealed on two timing
schedules, so that at every length one prefix determines the other. This is the form used in
Section~\ref{sec:fo}.
Proposition~\ref{prop:nojoin} shows that for stage masks the natural candidate for a join is not a join.

\section{Independent sources: the deficiencies do not add}\label{sec:source}

Theorem~\ref{thm:dich} compares an overwritten real with a random real. Now let both reals be overwritten: let
$R\oplus Y$ be random and let $F,G$ be computable sets with counting functions $f,g$. As in
Lemma~\ref{lem:over}(b), $K(R_F\rest n,Y_G\rest n)\ge 2n-f(n)-g(n)-O(1)$. So Theorem~\ref{thm:basic}, which for a
common upper bound $w$ gives $m-K(m)\le f(N)+g(N)+O(1)$ at $N=w_m$, shows that $R_F\incomp Y_G$ when
$f(n)+g(n)\le\gamma\log n$ for some $\gamma<1$. In the other direction Theorem~\ref{thm:dich}(i) shows that $R_F$
and $Y_G$ are compatible when $\sum_n2^{-f(n)}<\infty$ or $\sum_n2^{-g(n)}<\infty$, because $Y_G\leK Y$. When
$f(n)$ and $g(n)$ are both about $\frac34\log n$, neither argument applies: the joint complexity of the two reals
is too low for the method of Sections~\ref{sec:dip} and~\ref{sec:gauge}. This section shows by a different argument that the two reals are
incompatible in this case too. What counts is each deficiency separately, not their sum.

$K^B$ denotes prefix-free complexity relative to the oracle $B$. For a string $p$ let $n_p$ be the integer with
binary expansion $1p$; so $2^{|p|}\le n_p<2^{|p|+1}$, and $p$ and $n_p$ are computable from each other. Two facts
are used.

\smallskip
\noindent(F1) \emph{If $R$ is random, then $K(R\rest n)\ge n+K^R(n)-O(1)$.} By the ample excess lemma of Miller and
Yu \cite{MY}, $\sum_n2^{n-K(R\rest n)}<\infty$. The function $n\mapsto 2^{n-K(R\rest n)}$ is lower semicomputable
relative to $R$, so it is bounded by a constant times $2^{-K^R(n)}$. (Compare \cite{Mi}.)

\smallskip
\begin{lemma}[(F2): Muchnik's theorem relative to an oracle]\label{lem:muchnik}
There is a constant $c_M$ such that for every oracle $B$ and every string $a$ of length $\mu\ge2$ there is a string
$p$ with
\[
 |p|\le K^B(a)+c_M\log\mu,\qquad K(p\mid a)\le c_M\log\mu,\qquad K^B(a\mid p)\le c_M\log\mu.
\]
\end{lemma}
\begin{proof}
We use the online matching graphs from the proof of Muchnik's theorem by Musatov, Romashchenko and Shen
\cite[Section~2]{MRS}: for all $n$ and $k$ there is a bipartite graph, computable uniformly in $n,k$, whose left
part is $\{0,1\}^n$, whose right part has at most $2^k n^{O(1)}$ vertices, in which every left vertex has at most
$n^{O(1)}$ neighbours, and in which any sequence of at most $2^k$ distinct left vertices, presented one at a time,
can be matched to pairwise distinct neighbours by a computable online procedure that assigns each vertex a
neighbour when it arrives.

Let $k=K^B(a)+1$. If $k>\mu$, take $p=a$: then $K^B(a)\ge\mu=|p|$, and the other two inequalities hold trivially.
So we may assume $k\le\mu$, which is the range of parameters in \cite{MRS}. Take the graph for $n=\mu$ and this
$k$. Using $B$, enumerate
without repetitions the strings $x$ of length $\mu$ with $K^B(x)<k$; there are fewer than $2^k$ of them, and $a$ is
among them. Feed them to the online matching procedure in the order of enumeration, and let $p$ be the right
vertex assigned to $a$, written in binary. Then $|p|\le k+O(\log\mu)$. Given $a$, the numbers $\mu$ and $k$, and the
position of $p$ in the list of neighbours of $a$, one obtains $p$ without using $B$; as $\mu$, $k$ and this position
can be given by a self-delimiting description of $O(\log\mu)$ bits, $K(p\mid a)=O(\log\mu)$. Given $p$, $\mu$ and $k$
and the oracle $B$, repeat the enumeration and the matching until some string is assigned the vertex $p$; since
the matching is injective, that string is $a$. So $K^B(a\mid p)=O(\log\mu)$. All constants are independent of $B$.
\end{proof}
Only the enumeration of the strings of low complexity uses the oracle; the graph and the matching procedure do
not. This is why the construction of \cite{MRS} relativizes in the form needed here.

\begin{lemma}\label{lem:source}
Let $R$ be random, let $x$ be a real and $\delta$ a function with $K(R\rest n)\le K(x\rest n)+\delta(n)$ for all $n$,
and suppose $x\leK u$. Then there is a constant $c$ such that for every $\mu$ and every string $p$ with
$n_p\ge\mu$,
\[
 K^R(p)\le\delta(n_p)+K(p\mid u\rest\mu)+K(u\rest\mu)-\mu+c\le\delta(n_p)+K(p\mid u\rest\mu)+K(\mu)+2c.
\]
\end{lemma}
\begin{proof}
Put $\nu=u\rest\mu$ and $n=n_p$. Concatenate a shortest program for $\nu$, a shortest program for $p$ given $\nu$,
and the $n-\mu$ bits $u(\mu)\cdots u(n-1)$; a machine that has computed $\nu$ and $p$ knows $\mu$ and $n$ and
hence how many further bits to read. So
\[
 K(u\rest n)\le K(\nu)+K(p\mid\nu)+n-\mu+O(1).
\]
On the other hand, by (F1), $K(x\rest n)\ge K(R\rest n)-\delta(n)\ge n+K^R(p)-\delta(n)-O(1)$. Now use
$K(x\rest n)\le K(u\rest n)+O(1)$ for the first inequality and $K(\nu)\le\mu+K(\mu)+O(1)$ for the second.
\end{proof}
The lemma says that at the length $n_p$ the real $u$ is compressible, given $n_p$, by as much as $u\rest\mu$
knows about $p$, while $x$ is compressible only by $\delta(n_p)$ plus what $R$ knows about $p$.

For $\delta=0$ and $p=u\rest\mu$ the lemma gives $K^R(u\rest\mu)\le K(\mu)+O(1)$ for every $K$-upper bound $u$ of
a random real $R$. This is known in a stronger form: by Miller \cite{Mi}, if $R$ is random and $R\leK u$ then
$K^R(\sigma)\le K^u(\sigma)+O(1)$ for all strings $\sigma$. The point of Lemma~\ref{lem:source} is that it also
applies to reals $x$ that are not random but close to a random real.

\begin{remark}\label{rem:self}
For an overwritten real the oracle $R$ can be replaced by the real itself. Let $x=R_F$ with $R$ random and $F$
computable. Then $2^{n-f(n)-K(x\rest n)}\le2^{n-K(R\rest n)+O(1)}$, so this function of $n$ is summable, and it is
lower semicomputable relative to $x$. Hence $K(x\rest n)\ge n-f(n)+K^x(n)-O(1)$, and the proof of
Lemma~\ref{lem:source} gives $K^x(p)\le f(n_p)+K(p\mid u\rest\mu)+K(\mu)+c$ for every $u\ge_Kx$. This is
stronger, because $K^R\le K^x+O(1)$. Accordingly Proposition~\ref{prop:source} below holds for $x=R_F$ with
$K^{x}(u\rest\mu)\le C\log\mu$: up to $O(\log\mu)$ bits, a prefix of an upper bound of $R_F$ contains no
information that $R_F$ does not have.
\end{remark}

\begin{proposition}[Upper bounds are simple relative to the source]\label{prop:source}
Let $R$ be random, $\alpha<1$, and let $x$ be a real with $K(R\rest n)\le K(x\rest n)+\alpha\log n+O(1)$ for all
$n\ge1$. If $x\leK u$, then there is $C$ such that $K^R(u\rest\mu)\le C\log\mu$ for all $\mu\ge2$.
\end{proposition}
\begin{proof}
Fix $\mu$, put $\nu=u\rest\mu$ and $k=K^R(\nu)$, and choose $p$ by (F2) for $a=\nu$ and $B=R$. Then
\[
 k\le K^R(p)+K^R(\nu\mid p)+O(1)\le K^R(p)+c_M\log\mu+O(1).
\]
If $n_p<\mu$ then $|p|<\log\mu$, so $K^R(p)\le 2\log\mu+O(1)$ and $k\le(c_M+2)\log\mu+O(1)$. Otherwise
Lemma~\ref{lem:source} applies with $\delta(n)=\alpha\log n+O(1)$. As $\log n_p\le|p|+1\le k+c_M\log\mu+1$ and
$K(\mu)\le 2\log\mu+O(1)$,
\[
 K^R(p)\le\alpha k+(2c_M+2)\log\mu+O(1).
\]
The two displayed inequalities give $(1-\alpha)k\le(3c_M+2)\log\mu+O(1)$.
\end{proof}

\begin{lemma}[Independent sources have no common information]\label{lem:common}
If $R\oplus Y$ is random, there is $d$ such that for all strings $\nu$
\[
 K(\nu)\le K^R(\nu)+K^Y(\nu)+9\log\bigl(K^R(\nu)+K^Y(\nu)+2\bigr)+d.
\]
\end{lemma}
\begin{proof}
For a string $\nu$ and a number $i$ let $A^i_\nu=\{R':K^{R'}(\nu)\le i\}$. These classes are $\Sigma^0_1$ uniformly
in $\nu$ and $i$, and $\sum_\nu\lambda(A^i_\nu)<2^{i+1}$, because for each $R'$ fewer than $2^{i+1}$ strings $\nu$
have $K^{R'}(\nu)\le i$. For numbers $i,j,s$ let $Q_{i,j,s}$ be the set of all $\nu$ with
$\lambda(A^i_\nu)>2^{-j-s}$. This set has fewer than $2^{i+j+s+1}$ elements and is computably enumerable uniformly
in $i,j,s$. So for $\nu\in Q_{i,j,s}$
\begin{equation}\label{eq:Q}
 K(\nu)\le i+j+s+2\log(i+2)+2\log(j+2)+2\log(s+2)+O(1).
\end{equation}
Let $\widehat A_\nu$ be the $\Sigma^0_1$ class obtained by enumerating $A^i_\nu$ as long as the measure does not
exceed $2^{-j-s}$; so $\lambda(\widehat A_\nu)\le2^{-j-s}$, and $\widehat A_\nu=A^i_\nu$ if $\nu\notin Q_{i,j,s}$.
Then $V_{i,j,s}=\bigcup_\nu\widehat A_\nu\times A^j_\nu$ is a $\Sigma^0_1$ class of pairs of measure less than
$2^{-j-s}\cdot2^{j+1}=2^{1-s}$. Put $s(i,j,r)=r+2\lceil\log(i+2)\rceil+2\lceil\log(j+2)\rceil+1$ and
$W_r=\bigcup_{i,j}V_{i,j,s(i,j,r)}$. Then $\lambda(W_r)\le2^{-r}\sum_{i,j}(i+2)^{-2}(j+2)^{-2}\le2^{-r}$, so
$(W_r)_r$ is a Martin-L\"of test for pairs, and there is $r$ with $(R,Y)\notin W_r$.

Now let $\nu$ be any string and put $i=K^R(\nu)$, $j=K^Y(\nu)$ and $s=s(i,j,r)$. Then $R\in A^i_\nu$ and
$Y\in A^j_\nu$, but $(R,Y)\notin\widehat A_\nu\times A^j_\nu$. So $\widehat A_\nu\ne A^i_\nu$, hence
$\nu\in Q_{i,j,s}$, and \eqref{eq:Q} gives the bound.
\end{proof}

\begin{theorem}\label{thm:two}
Let $R\oplus Y$ be random, let $\alpha,\beta<1$, and let $x,y$ be reals such that for all $n\ge1$
\[
 K(R\rest n)\le K(x\rest n)+\alpha\log n+O(1)\qquad\text{and}\qquad K(Y\rest n)\le K(y\rest n)+\beta\log n+O(1).
\]
Then $x\incomp y$.
\end{theorem}
\begin{proof}
Suppose $x,y\leK u$. By Proposition~\ref{prop:source}, applied to $(R,x)$ and to $(Y,y)$, there is $C$ with
$K^R(u\rest\mu)\le C\log\mu$ and $K^Y(u\rest\mu)\le C\log\mu$ for all $\mu\ge2$. By Lemma~\ref{lem:common},
$K(u\rest\mu)\le 2C\log\mu+9\log(2C\log\mu+2)+d$. But
\[
 K(u\rest\mu)\ge K(x\rest\mu)-O(1)\ge K(R\rest\mu)-\alpha\log\mu-O(1)\ge\mu-\log\mu-O(1).
\]
This is a contradiction for large $\mu$.
\end{proof}

\begin{corollary}\label{cor:two}
Let $R\oplus Y$ be random and let $F,G$ be computable sets whose counting functions satisfy
$f(n)\le\alpha\log n+O(1)$ and $g(n)\le\beta\log n+O(1)$ with $\alpha,\beta<1$. Then $R_F\incomp Y_G$.
\end{corollary}
\begin{proof}
$K(R\rest n)\le K(R_F\rest n)+f(n)+O(1)$, as in Lemma~\ref{lem:over}(b), and similarly for $Y$.
\end{proof}

The proofs of Lemma~\ref{lem:source} and Proposition~\ref{prop:source} use the hypothesis on $x$ only at the
lengths $n_p$ that occur. This gives a local form, which improves Lemma~\ref{lem:blocks}.

\begin{lemma}[Short common blocks]\label{lem:short}
There is a constant $C_0$ with the following property, in the notation of Section~\ref{sec:masks}. Let $x,y$ be
reals, and suppose that for infinitely many $t$, with $N=N_{t-1}$, we have $x(k)=\rho^1(k)$ and $y(k)=\rho^2(k)$
for all $k$ with $N\le k<2^{N+\lceil C_0\log N\rceil}$. Then $x\incomp y$.
\end{lemma}
\begin{proof}
Let $A=2c_M+3$ and $C_0=A+c_M+1$, with $c_M$ from (F2). Suppose $x,y\leK u$, fix such a $t$ with $N$ large, and put
$W=2^{N+\lceil C_0\log N\rceil}$; then $W\le N_t$. For $N\le n<W$ the string $\rho^1\rest n$ is computable from
$x\rest n$ and $\rho^1\rest N$, and $N$ is computable from $n$. So $K(x\rest n)\ge n+K^{\rho^1}(n)-N-O(1)$ by (F1),
and the proof of Lemma~\ref{lem:source} gives a constant $c_1$, which does not depend on $t$, such that
\begin{equation}\label{eq:local}
 K^{\rho^1}(p)\le N+K(p\mid u\rest\mu)+K(\mu)+c_1\qquad\text{whenever }\max\{\mu,N\}\le n_p<W.
\end{equation}
Let $k(\mu)=K^{\rho^1}(u\rest\mu)$. There is a constant $c_0$ with $k(\mu+1)\le k(\mu)+c_0$ for all $\mu$.

\emph{Claim.} If $N$ is large, then $k(\mu)\le N+A\log N$ for $1\le\mu\le 4N$.

This is proved by induction on $\mu$. Let $2\le\mu\le4N$ and $k(\mu-1)\le N+A\log N$, so that
$k(\mu)\le N+A\log N+c_0$. Choose $p$ by (F2) for $a=u\rest\mu$ and $B=\rho^1$. Then
$|p|\le k(\mu)+c_M\log\mu$ and
\begin{equation}\label{eq:kmu}
 k(\mu)\le K^{\rho^1}(p)+c_M\log\mu+O(1).
\end{equation}
If $n_p<\max\{\mu,N\}$, then $|p|\le\log(4N)$, so $K^{\rho^1}(p)\le2\log(4N)+O(1)$ and \eqref{eq:kmu} gives
$k(\mu)\le(c_M+2)\log(4N)+O(1)\le N$. Otherwise
$n_p<2^{|p|+1}\le2^{N+A\log N+c_0+c_M\log(4N)+1}\le W$, since $N$ is large. So \eqref{eq:local} applies and gives
$K^{\rho^1}(p)\le N+c_M\log\mu+K(\mu)+c_1\le N+(c_M+2)\log(4N)+c_1+O(1)$. With \eqref{eq:kmu},
$k(\mu)\le N+(2c_M+2)\log(4N)+c_1+O(1)\le N+A\log N$. This proves the claim.

The same holds for $\rho^2$ and $y$. So $K^{\rho^1}(u\rest4N)$ and $K^{\rho^2}(u\rest4N)$ are at most $N+A\log N$,
and by Lemma~\ref{lem:common}
\[
 K(u\rest4N)\le2N+2A\log N+9\log(2N+2A\log N+2)+d.
\]
On the other hand $N\le4N<W$, so $K(u\rest4N)\ge K(x\rest4N)-O(1)\ge3N-O(1)$. For large $N$ this is a
contradiction.
\end{proof}

\begin{remark}\label{rem:two}
(1) Theorem~\ref{thm:two} does not require $x$ and $y$ to be computable from $R$ and $Y$. In
Corollary~\ref{cor:two}, for example, the positions in $F$ and $G$ may be overwritten by the bits of arbitrary
reals instead of zeros.

(2) Lemma~\ref{lem:short} is close to optimal. A real $x$ that copies $\rho^1$ only on
$[N_{t-1},2^{N_{t-1}/2})$ in each stage, and is $0$ elsewhere, satisfies $K(x\rest n)\le n+O(1)$ and so is
$K$-below every random real. So the block length that decides between the two cases is $2^{cN}$ with
$\frac12\le c\le1+o(1)$.
\end{remark}

Corollary~\ref{cor:two} and Theorem~\ref{thm:dich}(i) leave a narrow boundary case.
\begin{question}\label{q:two}
Let $R\oplus Y$ be random, and let $F,G$ be computable sets such that $\sum_n2^{-f(n)}$ and $\sum_n2^{-g(n)}$ both
diverge. Are $R_F$ and $Y_G$ incompatible?
\end{question}
If the answer is yes, the picture for this family is complete: $R_F$ and $Y_G$ are then compatible iff one of
the two reals is $K$-below every random real. Theorem~\ref{thm:dich}(ii) is the case $G=\emptyset$.

\smallskip
\noindent\textbf{Overwritten reals from one source.}
Question~\ref{q4} asks when $R_G\leK R_{G'}$ for one random real $R$. For sparse $G$ the answer is inclusion,
provided $R$ is \emph{weakly low for $K$} in the sense of Miller \cite{Mi}: there is $d$ such that
$K^R(j)\ge K(j)-d$ for infinitely many $j$. By \cite{Mi} this holds iff $\Omega$ is random relative to $R$; in
particular it holds for every $2$-random real.

\begin{proposition}\label{prop:order}
Let $R$ be random and weakly low for $K$. Let $F$ be a computable set such that $f(n)\le\log n+O(1)$ and
$f(2n)-f(n)\le O(1)$ for all $n\ge1$; for example $F=\{2^i:i\in\w\}$. Let $G\subseteq F$ and $G'$ be computable sets,
with counting functions $g,g'$. Then $R_G\leK R_{G'}$ iff $G'\subseteq^*G$.
\end{proposition}
\begin{proof}
If $G'\subseteq^*G$, then $R_G\rest n$ is computable from $R_{G'}\rest n$ with finitely many bits hard-wired.
Conversely, suppose $R_G\leK R_{G'}$ and let $D=G'\setminus G$ and $E=G\setminus G'$. We show that $D$ is finite.
Fix $d$ such that $J=\{j:K^R(j)\ge K(j)-d\}$ is infinite; in particular $J$ is unbounded, which is what is used
below.

\emph{Case 1: $E$ is finite.} Changing $G'$ at finitely many positions we may assume $G\subseteq G'$. By (F1),
$K(R_G\rest n)\ge K(R\rest n)-g(n)-O(1)\ge n+K^R(n)-g(n)-O(1)$, and $K(R_{G'}\rest n)\le K(n)+n-g'(n)+O(1)$ by
Lemma~\ref{lem:over}(c). So $|D\cap[0,n)|=g'(n)-g(n)\le K(n)-K^R(n)+O(1)$ for all $n$, and this is at most
$d+O(1)$ for $n\in J$.

\emph{Case 2: $E$ is infinite.} Let $e_1<e_2<\cdots$ be the elements of $E$. For $j\ge1$ let
$\sigma_j=R(e_1)\cdots R(e_j)$, and let $n_j$ be the integer with binary expansion $1\,\sigma_j\,1\,0^{r}$, where
$r$ is least such that this integer exceeds $e_j$. Then $j$ and $\sigma_j$ are computable from $n_j$, and $n_j$ is
computable from $\sigma_j$. The hypotheses on $F$ give a constant $b_1$ such that at most $b_1$ elements of $F$
lie between $e_j$ and $n_j$: if $r>0$ then $n_j\le2e_j$, and if $r=0$ then $n_j<2^{j+2}$, while
$j\le f(e_j+1)\le\log(e_j+1)+O(1)$. So $|E\cap[0,n_j)|\le j+b_1$.

The positions $e_i$ are not in $G'$. So $R_{G'}\rest n_j$ can be described by a shortest program for $n_j$,
followed by the bits of $R$ at the positions below $n_j$ that are neither in $G'$ nor among $e_1,\dots,e_j$:
\[
 K(R_{G'}\rest n_j)\le K(n_j)+n_j-g'(n_j)-j+O(1).
\]
The positions $e_i$ are in $G$. So from $R_G\rest n_j$, whose length gives $n_j$ and hence $\sigma_j$, and the
bits of $R$ at the other $g(n_j)-j$ positions of $G$ below $n_j$ one obtains $R\rest n_j$; by (F1)
\[
 K(R_G\rest n_j)\ge K(R\rest n_j)-g(n_j)+j-O(1)\ge n_j+K^R(n_j)-g(n_j)+j-O(1).
\]
As $K(R_G\rest n_j)\le K(R_{G'}\rest n_j)+O(1)$, we get $2j+g'(n_j)-g(n_j)\le K(n_j)-K^R(n_j)+O(1)$. Now
$K(n_j)\le j+K(j)+O(1)$ and $K^R(n_j)\ge K^R(j)-O(1)$, and
$g'(n_j)-g(n_j)=|D\cap[0,n_j)|-|E\cap[0,n_j)|\ge|D\cap[0,n_j)|-j-b_1$. Hence
\begin{equation}\label{eq:Dsmall}
 |D\cap[0,n_j)|\le K(j)-K^R(j)+O(1)\qquad\text{for all }j,
\end{equation}
and the right-hand side is at most $d+O(1)$ for $j\in J$.

In both cases $D$ has boundedly many elements below infinitely many lengths, so $D$ is finite.\qedhere
\end{proof}
For an arbitrary random $R$ the proof still gives \eqref{eq:Dsmall}, so that $G'\setminus G$ has at most
$2\log j+O(1)$ elements below $n_j$. We do not know whether the proposition holds for $R=\Omega$.

\begin{corollary}\label{cor:interval}
Let $R$ and $F$ be as in Proposition~\ref{prop:order}, with $F$ infinite. Then $G\mapsto R_G$ induces an
order-reversing embedding of the Boolean algebra of computable subsets of $F$ modulo finite sets into the
interval $[R_F,R]$ of $\mathcal D_K$. If moreover $f(n)\le\alpha\log n+O(1)$ with $\alpha<1$, for example if
$F=\{4^i:i\in\w\}$, then every $u\ge_KR_F$ satisfies $K^{R_F}(u\rest\mu)\le C\log\mu$ for some $C$ and all
$\mu\ge2$.
\end{corollary}
\begin{proof}
The first part is Proposition~\ref{prop:order}. The second is Proposition~\ref{prop:source} in the form of
Remark~\ref{rem:self}.
\end{proof}

So this is an interval of $\mathcal D_K$, just below a random degree, which contains a copy of the countable
atomless Boolean algebra and for which one lemma says something about \emph{every} element, and even about
every upper bound of its least element. Apart from the lower cones of Section~\ref{sec:sparse}, it is the only
definable set in this paper about all of whose elements something is proved. The statement is not precise enough
for a coding: a set $G\subseteq F$ has about $\log\mu$ elements below $\mu$, so the Boolean algebra lives inside
the error term $C\log\mu$.

\begin{question}\label{q:lub}
In the situation of Proposition~\ref{prop:order}, is $R_{G\cap G'}$ the least upper bound of $R_G$ and $R_{G'}$ in
$\mathcal D_K$, and $R_{G\cup G'}$ their greatest lower bound, for computable $G,G'\subseteq F$? Is every
complemented element of the interval $[R_F,R]$ the degree of some $R_G$?
\end{question}
For upper bounds of the form $R_H$ with $H$ computable the first part follows from
Proposition~\ref{prop:order}, which gives $H\subseteq^*G\cap G'$. The complexity function of $R_{G\cap G'}$ is not
the pointwise maximum of those of $R_G$ and $R_{G'}$: at most lengths it exceeds both by about
$\min\{|G\setminus G'|,|G'\setminus G|\}$ counted below $n$. So Lemma~\ref{lem:maxmin} does not apply, and a
positive answer would give least upper bounds of a different kind from those of Section~\ref{sec:timing}.

\section{Toward the first-order theory}\label{sec:fo}

\noindent\textbf{Status.}
For the Turing degrees the question is settled: the theory is undecidable (Lachlan \cite{Lach}) and is as
complicated as second-order arithmetic (Simpson \cite{Simp}). For $\mathcal D_K$ we know of no decidability or
undecidability result. The nearest published result is the undecidability of the Solovay degrees of c.e.\ reals
\cite{DHL}. $\mathcal D_K$ is not an upper semilattice, by \cite{MY}.

\smallskip
\noindent\textbf{What is proved here.}
The existential theory of $(\mathcal D_K,\le,\incomp)$ is decidable (Corollary~\ref{cor:fin}); this is the first
rung of the usual ladder. Undecidability of the full theory needs a \emph{definable} family, which no embedding
theorem provides. Section~\ref{sec:timing} supplies one ingredient of such a family: a Boolean algebra whose joins,
meets and complements are those of $\mathcal D_K$. We expect the theory to be undecidable.
Each of the two routes below, and the third one in Section~\ref{sec:profiles}, is missing a lemma that controls
\emph{all} the degrees in some definable set.

\smallskip
\noindent\textbf{Route A: initial segments.}
For the Turing degrees, and for several strong reducibilities, undecidability is proved by realizing every finite
distributive lattice as an initial segment $[0,a]$: the class of finite distributive lattices has a hereditarily
undecidable theory \cite{ELTT}, and $[0,a]$ is definable from $a$.
For $\mathcal D_K$ already the first case is open. Whether there is a minimal $K$-degree is a question from
Barmpalias's survey \cite{Bar}, and \cite{BL2} records it as unsolved.
Section~\ref{sec:sparse} at least turns it, for sparse codes, into a question about a strong reducibility, for
which tree constructions exist: Raichev and Stephan \cite{RS} construct a real of minimal $rK$-degree.

\begin{question}\label{q:min}
Is there a real $V$ such that $S_V$ has minimal $K$-degree? By Proposition~\ref{prop:sparse} this asks for a
non-$K$-trivial $V$ such that every real $y$ with $\sup_tK(y\rest4^t\mid(V\rest t)^*)<\infty$ is $K$-trivial or
satisfies $S_V\leK y$.
\end{question}
\begin{question}\label{q:lat}
Is every finite distributive lattice isomorphic to an initial segment of $\mathcal D_K$? A positive answer makes
the theory of $\mathcal D_K$ undecidable.
\end{question}
The difficulty to expect: below $S_V$ there are, besides the codes $S_W$, reals that reveal the bits of $V$ on a
delayed schedule, and for a general $V$ these form many degrees. $V$ must be built so that every delay either
loses nothing or loses everything; a real of minimal $rK$-degree has this flavour, since it must be predictable
from its own shorter prefixes.

The following elementary observation, which may well be known, shows how restrictive this is. For a computable
strictly increasing sequence $p=(p_k)_{k\in\w}$ with $p_k\ge k$ and a real $y$, let $S^p_y$ be the real with
$S^p_y(p_k)=y(k)$ for all $k$ and $S^p_y(n)=0$ for all other $n$. So $S_V=S^p_V$ for $p_k=4^k$.

\begin{proposition}[Dilution]\label{prop:dilute}
\begin{enumerate}[label=\textup{(\alph*)}]
\item $S^p_y\leK y$, and if $S^p_y$ is $K$-trivial then so is $y$. So below every nonzero $K$-degree lie the
nonzero degrees of its dilutions.
\item If $y\leK S^p_y$ for every such $p$, in particular if the $K$-degree of $y$ is minimal, then
$K(y\rest n)\le K(n)+h(n)+O(1)$ for every computable, nondecreasing, unbounded $h:\w\to\w$.
\end{enumerate}
\end{proposition}
\begin{proof}
(a) $S^p_y\rest n$ is computable from $y\rest n$, because $p_k\ge k$. And $y\rest k$ is computable from
$S^p_y\rest p_k$; so if $S^p_y$ is $K$-trivial, then $K(y\rest k)\le K(p_k)+O(1)\le K(k)+O(1)$.

(b) Let $t(n)=|\{k:p_k<n\}|$. Then $S^p_y\rest n$ is computable from $n$ and $y\rest t(n)$, so
\[
 K(y\rest n)\le K(S^p_y\rest n)+O(1)\le K(n)+K(y\rest t(n))+O(1)\le K(n)+3t(n)+O(1).
\]
Given $h$, let $p_k=k+\min\{n:h(n)\ge3(k+1)\}$. Then $3t(n)\le h(n)$.
\end{proof}
So a real of minimal $K$-degree, if there is one, lies in the narrow zone between trivial and nontrivial initial
segment complexity that is studied in \cite{BaBa}: its complexity exceeds $K(n)$ by an unbounded amount that
grows more slowly than every computable unbounded function. In Question~\ref{q:min} the real $V$ must have the
same property. Every real outside this zone has a nonzero degree strictly below its own.

\smallskip
\noindent\textbf{Route B: Nies's theorem.}
A \emph{$\Sigma^0_k$-Boolean algebra} is a Boolean algebra presented on the natural numbers with computable
operations and a $\Sigma^0_k$ ordering. It is \emph{effectively dense} if there is a function
$F\le_T\emptyset^{(k-1)}$ with $F(x)\le x$ for all $x$ and $0<F(x)<x$ for all $x\ne0$. Nies \cite{Nies97} proved that
for such an algebra $\mathcal B$ the lattice $\mathcal I(\mathcal B)$ of its $\Sigma^0_k$-ideals has a hereditarily
undecidable theory; it even interprets true arithmetic \cite{Nies00}. The theorem relativizes to any oracle.
Nies allows the Boolean operations to be computable in the relevant jump; our definition, which requires
computable operations, is a special case, so the sufficient condition below is stronger than necessary.

This is the engine of the undecidability proofs for degree structures on computable sets \cite{DN,Nies00}. As we
read \cite[Section~3]{Nies00}, the argument there has two halves. Let $A$ be a super sparse set and $a$ its
polynomial-time degree. First, the complemented elements of $[0,a]$ form a Boolean algebra $\mathcal B(a)$, whose
elements are the degrees of the sets $A\cap U$ with $U$ polynomial-time computable, and which is effectively
dense. Secondly, if $A$ is strongly super sparse, the sets $\{x\in\mathcal B(a):x\le c\}$ with $c\le a$ are exactly
the $\Sigma^0_2$-ideals of $\mathcal B(a)$. So $\mathcal I(\mathcal B(a))$ is definable in the degrees below $a$.

For $\mathcal D_K$ the two halves take the following form, with $\mathrm{Comp}[a,b]$ as defined before
Corollary~\ref{cor:comp}.

\begin{proposition}\label{prop:nies}
Suppose there are reals $a\leK b$ and $q$, an oracle $X$ and a number $k$ such that
\begin{enumerate}[label=\textup{(\alph*)}]
\item $\mathcal B=\mathrm{Comp}[a,b]$, with the order of $\mathcal D_K$, is a Boolean algebra which has a presentation
as an effectively dense $\Sigma^0_k(X)$-Boolean algebra;
\item the sets $I_c=\{y\in\mathcal B: y\le c\}$, for those $c\le q$ for which $I_c$ is an ideal of $\mathcal B$, are
exactly the $\Sigma^0_k(X)$-ideals of this presentation.
\end{enumerate}
Then the theory of $\mathcal D_K$ is undecidable. The same holds with an initial segment $[0,p]$ in place of
$\mathrm{Comp}[a,b]$.
\end{proposition}
\begin{proof}
$\mathcal B$ is definable from $a$ and $b$, and its join is definable from the order restricted to $\mathcal B$. So
``$a\le c\le q$ and $I_c$ is closed under the join of $\mathcal B$'' is a first-order condition on $c$ with
parameters $a,b,q$. On the set of such $c$ the definable preorder
``$\forall y\in\mathcal B\,(y\le c\to y\le c')$'' has a quotient isomorphic to $\mathcal I(\mathcal B)$. So this
lattice is definable with parameters in $\mathcal D_K$. By Nies's theorem relative to $X$ its theory is hereditarily
undecidable. The parameters are eliminated as in the proof of Proposition~\ref{prop:finite}: the sentences true in
all structures defined by parameters that satisfy the defining conditions form a deductively closed subtheory of a
hereditarily undecidable theory, so they cannot form a decidable set, and they are computable from the theory of
$\mathcal D_K$.
\end{proof}

\noindent\emph{What is known.} Take $a=z_\emptyset$ and $b=z_\w$. By Corollary~\ref{cor:comp} and
Remark~\ref{rem:timing}, $\mathrm{Comp}[a,b]$ contains the Boolean algebra $\mathcal B_0=\{z_S:S\text{ computable}\}$,
the operations of $\mathcal B_0$ are those of $\mathcal D_K$, and every ideal of $\mathcal B_0$ has the form
$\{y\in\mathcal B_0:y\le c\}$ for some degree $c$. This is the counterpart of the family of degrees of $A\cap U$.
Moreover $b$ vanishes on $[2N_{t-1}+1,N_t)$ for every $t$, so only countably many reals are $K$-below $b$
(Corollary~\ref{cor:countable}), and $\mathrm{Comp}[a,b]$ is countable.

\smallskip
\noindent\emph{What is missing.} Both hypotheses of Proposition~\ref{prop:nies} concern all of $\mathrm{Comp}[a,b]$,
and we can prove neither. In the polynomial-time proof the first half comes almost for free: a time-bounded
reduction to a super sparse set can compute the earlier part of the set by brute force, the degrees below $a$
then form a distributive lattice \cite{AS,DN}, and in a distributive lattice the complemented elements form a
Boolean algebra. Complexity has no time bound, earlier information never becomes free, and lower cones in
$\mathcal D_K$ are not lattices of this kind; Proposition~\ref{prop:nojoin} is an instance. The interval $[a,b]$
contains, besides the timing masks, reals that show part of a block early and the rest late, or that show a
block early after transforming it. One has to prove that none of these has a complement, and to find out which
sets $S$ occur.

\begin{question}\label{q:comp}
Let $a=z_\emptyset$ and $b=z_\w$. Is every element of $\mathrm{Comp}[a,b]$ the degree of a real $y$ such that
$K(y\rest n)=K(n,\rho\rest m_S(n))\pm O(1)$ for some set $S$? For which sets $S$ is there such a real?
\end{question}
All computable and all $K$-trivial sets $S$ occur, and so do some sets computed by $\rho$ itself, for instance
$S=\{t:\rho(\lfloor t/2\rfloor)=1\}$, since then Lemma~\ref{lem:level}(b) still holds. So the algebra in (a), if it
is one, is larger than $\mathcal B_0$ and depends on $\rho$; the level $k$ and the oracle $X$ have to be read off
from the answer.

\begin{question}\label{q:exact}
For a suitable $\rho$, is there a real $q$ such that the ideals $I_c$ with $c\le q$ are exactly the ideals of one
arithmetical level of $\mathrm{Comp}[z_\emptyset,z_\w]$?
\end{question}
With arbitrary degrees as names one gets every ideal (Remark~\ref{rem:timing}), and then only the theory of
a countable Boolean algebra with quantification over ideals, which is decidable by Rabin's theorem \cite{Rabin}.
The restriction to a countable definable set of names, such as the reals below a real with gaps, has to make the
match exact. We have not found the right $q$.

\begin{question}\label{q:def}
Is there a real $p$ such that the $K$-degrees below $p$ form an infinite Boolean algebra?
\end{question}
This is the initial-segment form of hypothesis (a). If the algebra has atoms, they are minimal $K$-degrees.

\begin{remark}\label{rem:closure}
For hypothesis (a) a classification of the elements is more than is needed. Suppose that every function obtained
from the complexity functions $n\mapsto K(x\rest n)$ of the elements $x$ of $\mathrm{Comp}[a,b]$ by finitely many
pointwise maxima and minima is, up to an additive constant, again the complexity function of a real. Then by
Lemma~\ref{lem:maxmin} these reals form a distributive lattice under $\le$ in which the generators have
complements, hence a Boolean algebra, and it equals $\mathrm{Comp}[a,b]$. So it would suffice to prove this closure
property. We do not see how to prove it without knowing the elements. Note that Lemma~\ref{lem:maxmin} gives only a
sufficient condition: we do not know whether every least upper bound or greatest lower bound in $\mathcal D_K$ is a
pointwise maximum or minimum (this is question (P1) in Section~\ref{sec:proposed}).
\end{remark}

\begin{remark}[Decidable fragments]\label{rem:fragments}
On the decidable side what can be proved at present is fragments.

(1) The existential theory of $(\mathcal D_K,\le,\incomp)$ is decidable (Corollary~\ref{cor:fin}(c)). As $\incomp$ is
defined by a universal formula, this is slightly more than the existential theory of the order.

(2) The next level contains open problems. ``There is a minimal degree'' is the $\exists\forall$ sentence
\[
 \exists x\,\exists w\,\forall y\,\forall z\;\bigl[\,x\not\le w\ \wedge\ \bigl(y\le x\to(x\le y\vee y\le z)\bigr)\bigr],
\]
and its truth value is not known (Question~\ref{q:min}). So a decision procedure for the two-quantifier
fragment, applied to this sentence, would settle the minimal degree problem; and one for the full theory would
also settle whether any two degrees have a greatest lower bound \cite[Question~9.1]{MN}.

(3) Between (1) and (2) lies the existential theory in the language with the two ternary relations ``$z$ is the
least upper bound of $x,y$'' and ``$z$ is the greatest lower bound of $x,y$''. By Theorem~\ref{thm:timing} every
finite distributive lattice embeds into $\mathcal D_K$ with least upper bounds and greatest lower bounds
preserved, because it is a sublattice of a finite Boolean algebra and hence of the Boolean algebra of
computable sets modulo finite sets. The first sentences of this fragment whose truth we cannot determine are
the following.
\end{remark}

\begin{question}\label{q:M3}
Do the five-element nondistributive lattices $N_5$ and $M_3$ embed into $\mathcal D_K$ with least upper bounds
and greatest lower bounds preserved?
\end{question}
Lemma~\ref{lem:maxmin} cannot produce such an embedding: pointwise maxima and minima of functions satisfy the
distributive laws, so a finite set of degrees in which all least upper bounds and greatest lower bounds are
given by Lemma~\ref{lem:maxmin} generates a distributive lattice. A candidate for $M_3$ of a different kind is the
following. Let $R$ and $F$ be as in Proposition~\ref{prop:order}, let $d_1<e_1<d_2<e_2<\cdots$ be the elements of
$F$, and put $F_1=\{e_i\}$ and $F_2=\{d_i\}$. Let $a=R_{F_1}$, $b=R_{F_2}$, and let $c$ agree with $R$ outside $F$,
with $c(e_i)=R(e_i)+R(d_i)\bmod2$ and $c(d_i)=0$. All three reals lie in the interval $[R_F,R]$, each of them
lacks one bit of $R$ for every $i$, and any two of them together determine $R$. If $R$ is the least upper bound
and $R_F$ the greatest lower bound of each of the three pairs, this is a copy of $M_3$. Whether that is so is a
special case of the problem behind Question~\ref{q:lub}: to show that every upper bound of two such reals is an
upper bound of $R$.

\smallskip
Within Route B, Question~\ref{q:comp} is the natural first step. It concerns one explicit interval; by Proposition~\ref{prop:rigid} every
real in it is computable, with bounded advice, from a shortest program for the prefix of $z_\w$ at the lengths
$2N_{t-1}+1$; and unlike Route A it does not contain the problem of minimal $K$-degrees.

\section{Finite structures and incompatibility profiles}\label{sec:profiles}

A third route is the following: interpret finite graphs instead of
an ideal lattice, and classify degrees only up to their incompatibility behaviour against a family of witnesses.
This avoids minimal degrees and Question~9.8. The section records what we can prove in this direction and where
it stops.

\begin{proposition}\label{prop:finite}
Let $\mathcal C$ be a class of finite structures in a finite relational language whose theory is hereditarily
undecidable, for example the finite graphs or the finite partial orders \cite{ELTT}. Suppose there are formulas
$\delta(x;\bar p)$, $\epsilon(x,y;\bar p)$ and $\theta_R(\bar x;\bar p)$, one for each relation symbol $R$, in the
language of partial orders, such that every structure in $\mathcal C$ is isomorphic to one of the following form,
for suitable parameters $\bar p$ in $\mathcal D_K$: the domain is $\{x:\delta(x;\bar p)\}$ modulo the equivalence
relation defined by $\epsilon$, and each $R$ is interpreted by $\theta_R$. Then the theory of $\mathcal D_K$ is
undecidable.
\end{proposition}
\begin{proof}
Let $\nu(\bar p)$ say that $\epsilon$ defines an equivalence relation on the nonempty set defined by $\delta$ and
that each $\theta_R$ respects it. For a sentence $\varphi$ in the language of $\mathcal C$ let $\varphi^*(\bar p)$ be
its translation. The set $T$ of all $\varphi$ with
$\mathcal D_K\models\forall\bar p\,(\nu(\bar p)\to\varphi^*(\bar p))$ is the theory of the class of all structures
obtained from parameters satisfying $\nu$. So $T$ is closed under consequence, and
$T\subseteq\mathrm{Th}(\mathcal C)$ because every member of $\mathcal C$ is among these structures. If the theory
of $\mathcal D_K$ were decidable, $T$ would be decidable, contradicting hereditary undecidability.
\end{proof}
So parameters for which the formulas define something unintended, even an infinite structure, do no harm. What
has to be controlled is every element of the defined domain for the intended parameters. The number of
parameters must not depend on the structure.

\smallskip
\noindent\textbf{Profiles.}
For a set $\mathcal W$ of degrees and a degree $x$ let $P_{\mathcal W}(x)=\{w\in\mathcal W:x\incomp w\}$ be the
\emph{profile} of $x$. Write $x\preceq_{\mathcal W}x'$ if $P_{\mathcal W}(x)\subseteq P_{\mathcal W}(x')$, and
$x\sim_{\mathcal W}x'$ if the profiles are equal. If $\mathcal W$ is definable with parameters, so are these two
relations, because $\incomp$ is definable from $\le$.

\begin{lemma}\label{lem:profile}
\begin{enumerate}[label=\textup{(\alph*)}]
\item If $x\le x'$ then $P_{\mathcal W}(x)\subseteq P_{\mathcal W}(x')$.
\item If $x$ and all elements of $\mathcal W$ lie below one degree, then $P_{\mathcal W}(x)=\emptyset$. In
particular the elements of one lower cone are pairwise compatible.
\item Every real $K$-below the real $z_\w$ of Section~\ref{sec:timing} is $K$-below every random real. So for
$x\leK z_\w$ the profile $P_{\mathcal W}(x)$ is empty whenever all elements of $\mathcal W$ have random upper
bounds, and, by Proposition~\ref{prop:dim} \textup{(}which uses \cite{Zhao}\textup{)}, for every $\mathcal W$.
\item Let $T$ be an infinite computable set of stages. If $\mathcal W$ contains the degree of $\rho^2$ masked to
$B$ for every computable $B\subseteq T$, then the reals $\rho^1$ masked to $A$, for computable $A\subseteq T$, have
different profiles unless the sets agree modulo finite sets.
\end{enumerate}
\end{lemma}
\begin{proof}
(a) is (C2), and (b) is clear. (d) If $A\setminus A'$ is infinite, take $B=A\setminus A'$ and apply
Theorem~\ref{thm:masks}(b) to $A,B$ and to $A',B$.

(c) Let $N_{t-1}<n\le N_t$. By Lemma~\ref{lem:level}(b), $K(z_\w\rest n)\le m_\w(n)+K(n)+O(1)$. Note that
$c_{t-1}\le N_{t-1}/2$ and $c_t\le 3N_{t-1}/2$. If $n\le 2N_{t-1}+1$ then $m_\w(n)\le c_{t-1}+n-N_{t-1}\le n-N_{t-1}/2$
and $K(n)\le 2\log N_{t-1}+O(1)$. If $n>2N_{t-1}+1$ then $m_\w(n)=c_t\le 3N_{t-1}/2$ and
$n-K(n)\ge 3N_{t-1}/2-O(1)$. In both cases $K(z_\w\rest n)\le n+O(1)$. So $K(x\rest n)\le n+O(1)$ for every
$x\leK z_\w$, and Proposition~\ref{prop:below} applies. The same estimates give $\limsup_nK(x\rest n)/n\le 3/4$.
\end{proof}

Part (c) says that the interval of Section~\ref{sec:timing} is invisible to incompatibility: timing masks can
serve neither as vertices nor as witnesses in a coding by profiles. Part (d) says that with computable supports
on both sides the quotient by profile equivalence is infinite.

For a fixed number $n$ of witnesses nothing has to be proved. Against witnesses $w_1,\dots,w_n$ the profile of a
degree is a subset of $\{1,\dots,n\}$; and if $B_1,\dots,B_n$ are disjoint infinite computable sets of stages and
$w_v$ is $\rho^2$ masked to $B_v$, then by Theorem~\ref{thm:masks}(b) all subsets occur as profiles of masks of
$\rho^1$. But an interpretation may use only a bounded number of parameters. So the set of witnesses has to be
definable from boundedly many parameters, and such a set, a lower cone for example, is infinite and contains
with each witness thinner and weaker ones. The next results show that thinner witnesses can be made harmless.
Weaker ones are discussed at the end of the section.

\smallskip
\noindent\textbf{Witnesses with arbitrary supports.}
For a set $D$ of stages write $x^i_D$ for $\rho^i$ masked to $D$: $x^i_D(k)=\rho^i(k)$ if $k\in P_t$ with $t\in D$,
and $x^i_D(k)=0$ otherwise.

\begin{proposition}\label{prop:arb}
Let $A$ be a computable set of stages and $D$ an arbitrary set of stages.
\begin{enumerate}[label=\textup{(\alph*)}]
\item If $A\cap D$ is infinite, then $x^1_A\incomp x^2_D$.
\item If $A\cap D$ is finite, then $x^1_A$ and $x^2_D$ are compatible.
\end{enumerate}
\end{proposition}
\begin{proof}
(a) is Lemma~\ref{lem:blocks}, applied to the stages in $A\cap D$. (b) Let $u$ copy $\rho^1$ on the stages in $A$
and $\rho^2$ on the stages in $D\setminus A$, and be $0$ elsewhere. Then $x^1_A\rest n$ is computable from
$u\rest n$ by erasing the stages outside $A$, and $x^2_D\rest n$ is computable from $u\rest n$, with finitely many
bits hard-wired, by erasing the stages in $A$ and restoring the finitely many stages in $A\cap D$.
\end{proof}

Recall that an infinite set $D$ is \emph{r-cohesive} if for every computable set $A$ one of $D\cap A$,
$D\setminus A$ is finite; every infinite set has an r-cohesive subset.

\begin{corollary}\label{cor:cohesive}
Let $D$ be r-cohesive. Then for every computable $A$ and every infinite $D'\subseteq D$:
$x^1_A\incomp x^2_{D'}$ iff $D\subseteq^*A$. So against the masks of computable sets all the reals $x^2_{D'}$ with
$D'\subseteq D$ infinite behave in the same way.
\end{corollary}

Corollary~\ref{cor:cohesive} controls the masks of computable sets. The next statement controls arbitrary reals,
by a soft argument.

\begin{proposition}[Zero--one law]\label{prop:zeroone}
Let $\mathcal Y$ be a countable set of reals and $B$ an infinite set of stages. Then there is an infinite
$D\subseteq B$ such that for every $y\in\mathcal Y$ one of the following holds:
\begin{enumerate}[label=\textup{(\roman*)}]
\item $y\incomp x^2_{D'}$ for every infinite $D'\subseteq D$;
\item $y$ and $x^2_{D'}$ are compatible for every infinite $D'\subseteq D$.
\end{enumerate}
\end{proposition}
\begin{proof}
For $y\in\mathcal Y$ let $\mathcal P_y$ be the set of infinite $D'\subseteq B$ such that $y$ and $x^2_{D'}$ are
compatible, that is, such that for some real $u$ and some $c$ we have $K(y\rest n)\le K(u\rest n)+c$ and
$K(x^2_{D'}\rest n)\le K(u\rest n)+c$ for all $n$. The condition on $(u,D')$ is arithmetical, so $\mathcal P_y$ is an
analytic set of infinite subsets of $B$. By Silver's theorem \cite{Silver} analytic sets are Ramsey: every infinite
subset of $B$ has an infinite subset $H$ such that either all infinite subsets of $H$ are in $\mathcal P_y$ or none
is. Let $\mathcal Y=\{y_0,y_1,\dots\}$, choose $B\supseteq H_0\supseteq H_1\supseteq\cdots$ with $H_k$ of this kind for
$y_k$, and let $D=\{d_0<d_1<\cdots\}$ with $d_k\in H_k$. If $D'\subseteq D$ is infinite then
$D'\setminus\{d_0,\dots,d_{k-1}\}$ is an infinite subset of $H_k$. Removing finitely many stages from $D'$ changes
$x^2_{D'}$ at finitely many positions and so does not change its degree. Hence $D'\in\mathcal P_{y_k}$ iff
$D\in\mathcal P_{y_k}$.
\end{proof}

\begin{theorem}\label{thm:Bn}
Let $A_0$ be an infinite computable set of stages with infinite complement, let $q=x^1_{A_0}$, let $n\ge1$, and
let $B_1,\dots,B_n$ be a partition of $A_0$ into infinite computable sets. Then there are infinite sets
$D_v\subseteq B_v$ such that, with
\[
 \mathcal W=\{x^2_{D'}: D'\subseteq D_v\text{ infinite},\ 1\le v\le n\}\qquad\text{and}\qquad
 \pi(y)=\{v: y\incomp x^2_{D_v}\},
\]
the following hold.
\begin{enumerate}[label=\textup{(\alph*)}]
\item For every real $y\leK q$, every $v$ and every infinite $D'\subseteq D_v$: $y\incomp x^2_{D'}$ iff $v\in\pi(y)$.
\item Every subset of $\{1,\dots,n\}$ is $\pi(y)$ for some $y\leK q$.
\item For $y,y'\leK q$: $y\preceq_{\mathcal W}y'$ iff $\pi(y)\subseteq\pi(y')$. So the lower cone of $q$ modulo
$\sim_{\mathcal W}$, ordered by $\preceq_{\mathcal W}$, is the Boolean algebra with $n$ atoms.
\end{enumerate}
\end{theorem}
\begin{proof}
Only countably many reals are $K$-below $q$, by Corollary~\ref{cor:countable}(a). Apply
Proposition~\ref{prop:zeroone} to this set of reals and to each $B_v$ in turn; this gives $D_v$ and (a).
For (b) let $V'\subseteq\{1,\dots,n\}$, let $A$ be the union of the $B_v$ with $v\in V'$, and $y=x^1_A$. Then $y\leK q$
by Theorem~\ref{thm:masks}(a), and by Proposition~\ref{prop:arb}, $y\incomp x^2_{D_v}$ iff $D_v\subseteq A$ iff
$v\in V'$. (c) follows from (a).
\end{proof}

What Theorem~\ref{thm:Bn} adds to the case of $n$ witnesses is that $\mathcal W$ is closed under thinning and
has the size of the continuum, and still every degree below $q$, whatever it is, has the profile of one of
$2^n$ designated masks. Nothing has to be known about the reals below $q$ except that there are countably many.
But $\mathcal W$ is not known to be definable from boundedly many parameters. The obvious definable set is the
lower cone of $r=x^2_{D_1\cup\dots\cup D_n}$, which contains each $x^2_{D_v}$. Its elements are not all of the form
$x^2_{D'}$, and the zero--one law does not apply to them as it stands: its proof needs witnesses that are given
uniformly from $D'$ and whose degree does not change when $D'$ is changed finitely. Already for the reals
$x^2_{D_v\cup D_{v'}}$, which are in the cone, we cannot show that a degree below $q$ which is compatible with
$x^2_{D_v}$ and with $x^2_{D_{v'}}$ is compatible with $x^2_{D_v\cup D_{v'}}$.

\smallskip
\noindent\textbf{A bound on the complexity of lower cones.}
\begin{proposition}\label{prop:turing}
Suppose $K(x\rest n)\le K(n)+O(1)$ for infinitely many $n$. Then every real $y\leK x$ is computable in
$x\oplus\emptyset'$.
\end{proposition}
\begin{proof}
Choose $a$ such that $L=\{n:K(x\rest n)\le K(n)+a\}$ is infinite, and let $K(y\rest n)\le K(x\rest n)+c$ for all $n$.
For $n\in L$ the string $y\rest n$ lies in $F_n=\{\sigma\in2^n:K(\sigma)\le K(n)+a+c\}$, which has at most
$2^{a+c+O(1)}$ elements. The set $L$ is computable in $x\oplus\emptyset'$ and the sets $F_n$ are uniformly computable
in $\emptyset'$. So the tree of all strings $\sigma$ with $\sigma\rest n\in F_n$ for all $n\in L$, $n\le|\sigma|$, is
computable in $x\oplus\emptyset'$ and has boundedly many nodes on the levels in $L$. It has finitely many infinite
paths, each is isolated, and an isolated path of a binary tree is computable in the tree.
\end{proof}
For such $x$ Barmpalias and Vlek \cite{BV} proved that all $K$-predecessors are $\Delta^0_2$ relative to $x$; the
counting argument above gives the oracle $x\oplus\emptyset'$. It has two uses here. Every real
below $z_\w$ is computable in $\rho\oplus\emptyset'$, and on indices of such reals $\leK$ is a $\Sigma^0_2$
relation relative to $\rho\oplus\emptyset'$; so $X=\rho\oplus\emptyset'$ is a natural candidate for the oracle in
Proposition~\ref{prop:nies}. This bounds only the complexity of the order; a presentation as in
Proposition~\ref{prop:nies} also requires an effective control of the domain $\mathrm{Comp}[a,b]$ and of its
operations, which we do not have. And every real below the real $q$ of Theorem~\ref{thm:Bn} is computable in
$\rho^1\oplus\emptyset'$.

\smallskip
\noindent\textbf{A target with two parameters.}
A Boolean algebra is not enough, because the theory of finite Boolean algebras is decidable; a binary relation
is needed. Incompatibility between two lower cones provides one. For degrees $q_1,q_2$ and $x\le q_1$, $w\le q_2$
put
\[
 P(x)=\{w'\le q_2:x\incomp w'\}\qquad\text{and}\qquad P'(w)=\{x'\le q_1:x'\incomp w\}.
\]
Call $x\le q_1$ an \emph{atom} of $(q_1,q_2)$ if $P(x)\ne\emptyset$ and there is no $x'\le q_1$ with
$\emptyset\ne P(x')\subsetneq P(x)$; atoms $w\le q_2$ are defined in the same way with $P'$. A degree below both
$q_1$ and $q_2$ is not an atom, by Lemma~\ref{lem:profile}(b). Let $\mathcal G(q_1,q_2)$ be the bipartite graph
whose left vertices are the atoms below $q_1$ modulo equality of $P$, whose right vertices are the atoms below
$q_2$ modulo equality of $P'$, and in which the class of $x$ is adjacent to the class of $w$ iff $x\incomp w$.
This does not depend on the representatives, since $x\incomp w$ iff $w\in P(x)$ iff $x\in P'(w)$.

\begin{proposition}\label{prop:target}
Suppose that for every finite graph $G$ in which every vertex has degree at least $2$ there are degrees
$q_1,q_2$ such that $\mathcal G(q_1,q_2)$ is isomorphic to the incidence graph of $G$ \textup{(}left vertices the
vertices of $G$, right vertices the edges of $G$, and a vertex adjacent to the edges that contain
it\textup{)}. Then the theory of $\mathcal D_K$ is undecidable.
\end{proposition}
\begin{proof}
$\mathcal G(q_1,q_2)$ is defined in $\mathcal D_K$ by fixed formulas with the parameters $q_1,q_2$. The theory of
finite graphs is hereditarily undecidable \cite{ELTT}. An arbitrary finite graph is definable with two parameters
in the graph obtained by adding two new adjacent vertices joined to all old ones, and in the new graph every
vertex has degree at least $2$. A graph is definable in its incidence graph: two vertices are adjacent iff they
are distinct and have a common neighbour. So the class of incidence graphs in the hypothesis has a hereditarily
undecidable theory, by the argument of Proposition~\ref{prop:finite}, and Proposition~\ref{prop:finite} applies.
\end{proof}

The restriction on $G$ is there because in the incidence graph of such a graph no two vertices on the same side
have comparable neighbourhoods. The intended $q_1,q_2$ for a bipartite graph $(V_1,V_2,E)$ of this kind are as
follows. Choose disjoint infinite sets $S_e$ of stages, one for each $e\in E$, all contained in one r-cohesive
set, with infinitely many stages in no $S_e$. Let $q_i$ copy $\rho^i$ on the stages in
$T_i=\bigcup_eS_e$, and let the last positions of stage $t-1$ in $q_i$ say whether $t\in T_i$ and, if so, which
vertex in $V_i$ is the endpoint of the edge $e$ with $t\in S_e$. Concretely, the label of stage $t$ is written, in a fixed self-delimiting code, on the last
$O(\log|V_i|)$ positions of $P_{t-1}$, after the copied block of stage $t-1$ (if any) and after a run of zeros of
length at least $N_{t-1}/2$, so that the lengths before the labels are announced in the sense of
Section~\ref{sec:sparse}. A decoder given $q_i\rest n$ determines the stage boundaries from $n$, reads every complete
label, and reads the copied blocks of the stages labelled active. For $A\subseteq V_i$ let $q_i^A$ be $q_i$ with the copy
of $\rho^i$ erased on all stages whose label is not in $A$; the labels are kept. Then $q_i^A\leK q_i$, because the
label of a stage is visible before the stage begins. If some $a\in A\subseteq V_1$ is adjacent to some $b\in B\subseteq V_2$, then $q_1^A\incomp q_2^B$ by
Lemma~\ref{lem:blocks}; otherwise the two reals are active on disjoint sets of stages and a real that carries
both sequences of labels and copies each source on its stages is a common upper bound. So the reals $q_1^A$ and
$q_2^B$ realize the incidence pattern, and the profiles of the $q_1^{\{a\}}$ are pairwise incomparable. Unlike
for the computable supports of Lemma~\ref{lem:profile}(d), the mask of a single $S_e$ need not lie in either
cone: a computable set contains almost all or almost none of the stages in $T_i$, the labels on one side do
not determine the edge, and if the $S_e$ are chosen so that $K(S_e\cap[0,t)\mid(q_i\rest N_{t-1})^*)$ is unbounded
on the stages $t\notin T_i$, then $x^i_{S_e}\not\leK q_i$ by Proposition~\ref{prop:rigid}.

\smallskip
\noindent\textbf{What is missing.}
For the hypothesis of Proposition~\ref{prop:target} one has to show that the $q_i^{\{a\}}$ are atoms and that
there are no other atoms up to equality of profiles. This is a statement about \emph{all} degrees below $q_1$
and below $q_2$, and we cannot prove it. What is known is that they are countable and computable in
$q_i\oplus\emptyset'$ (Proposition~\ref{prop:turing}) and rigid at announced lengths (Proposition~\ref{prop:rigid}),
and that weakened elements behave like full ones away from a narrow boundary (Theorem~\ref{thm:two},
Lemma~\ref{lem:short} and Remark~\ref{rem:two}(2)); the boundary case is Question~\ref{q:two}. A two-sided version of
the zero--one law (Proposition~\ref{prop:zeroone}) would be the natural tool, but names for the elements of a cone are
not invariant under finite changes of the support.

There is also a general limitation, which any proof has to get around. The masks of computable sets behave like
the elements of the Boolean algebra $\mathbb B$ of computable sets modulo finite sets, and a parameter $x^i_D$
acts on them only through the ideals $\{A:A\subseteq^*D\}$ (Remark~\ref{rem:ideal}) and
$\{A:A\cap D\text{ is finite}\}$ (Proposition~\ref{prop:arb}). The theory of $\mathbb B$ with quantification over
ideals is decidable \cite{Rabin}. So a construction in which vertices, witnesses and parameters enter only through
this calculus cannot yield a hereditarily undecidable class of finite structures. The labels above, which are
not computable from the stage number, are meant to supply what the calculus lacks.

\begin{question}\label{q:cone}
Let $q_1,q_2$ be as described after Proposition~\ref{prop:target}. Are there only finitely many profiles $P(x)$,
$x\le q_1$? Are the classes of the reals $q_1^{\{a\}}$ exactly the atoms?
\end{question}

Of the open problems in Sections~\ref{sec:fo}--\ref{sec:profiles}, Question~\ref{q:two} is the most concrete: it
concerns two explicit reals, and its answer decides whether overwritten blocks can produce new profiles.

\section{Open questions, and comparison with other degree structures}\label{sec:open}

This section collects open questions about $\mathcal D_K$: first those in the literature, with their status as far as we could
determine it in November 2026, and then new questions. \emph{Answered} means
answered in print; \emph{open} means that we found no answer in the literature.

\subsection{Questions from the literature}
\noindent\emph{Miller and Nies \cite[Section~9]{MN}.}
\begin{enumerate}[label=\textup{(L\arabic*)},leftmargin=3em]
\item \cite[Q.~9.1]{MN} Do meets always exist in the vL-, $K$- and $C$-degrees? Is there a non-zero branching
degree? \emph{Answered for the $K$- and $C$-degrees by Merkle and Stephan} \cite{MS}: some pairs have no
meet, and there are non-zero branching degrees; Corollary~\ref{cor:branch} gives explicit examples with both
meet and join. \emph{Open for the vL-degrees.}
\item \cite[Q.~9.2]{MN} For $n\ge2$, are the weakly $n$-random, Schnorr $n$-random and computably $n$-random reals
closed upward in the $K$-degrees? \emph{Open}; the case of weak $2$-randomness, a question of
Miller and Yu, is listed as open in \cite[Open Question~11]{Blog17}.
\item \cite[Q.~9.3]{MN} Is there a natural characterization of computable $2$-randomness by initial segment
complexity? \emph{Open.}
\item \cite[Q.~9.6]{MN} Are there comparable $K$-degrees of $2$-random reals? \emph{Answered positively in the
preprint} \cite[Theorem~1.4]{Liu}.
\item \cite[Q.~9.7]{MN} Are there maximal $K$-degrees, or maximal random $K$-degrees? \emph{Answered negatively
by Liu} \cite{Liu}; this also refutes the conjecture in \cite{Mi} that $2$-random reals have maximal $K$-degree.
\item \cite[Q.~9.8]{MN} Is every real $K$-below a random real? \emph{Open}; see Sections~\ref{sec:over}--\ref{sec:fo}
of this paper and \cite{Zhao}. The analogue for $rK$ fails \cite{BL}.
\item \cite[Q.~9.9]{MN} For random $X,Y$, does $X\equiv_KY$ imply $X\equiv_TY$? \emph{Answered negatively by
Merkle and Stephan} \cite[Theorem~2 and the discussion following it]{MS}.
\item \cite[Q.~9.10]{MN} Are there random $A<_CB$? \emph{Answered positively by Liu} \cite{Liu}.
\item \cite[Q.~9.11]{MN} Does $X\le_CY$ imply $X\leK Y$? Do $\le_K$ and $\le_C$ differ on random reals?
\emph{Open.}
\end{enumerate}

\noindent\emph{Barmpalias's survey \cite[Section~3 and later]{Bar}.}
\begin{enumerate}[label=\textup{(L\arabic*)},leftmargin=3em,resume]
\item \cite[Q.~1]{Bar} Are the infinitely often $K$-trivial and the infinitely often $C$-trivial reals different
classes? \emph{Answered positively by Bauwens} \cite{Bau}: there is a real $X$ with
$\liminf_n(C(X\rest n)-C(n))=\infty$ and $\liminf_n(K(X\rest n)-K(n))<\infty$.
\item \cite[Q.~4]{Bar} Are there random $X<_KY$ with $\liminf_n|K(Y\rest n)-K(X\rest n)|<\infty$? \emph{Answered positively by Liu} \cite[Theorem~1.4]{Liu}, who also answers the related question of Miller
and Yu whether $x<_Ky$ forces $K(y\rest n)-K(x\rest n)\to\infty$.
\item \cite[Q.~5]{Bar} Are there non-$K$-trivial $X,Y$ and $c$ with
$\min\{K(X\rest n),K(Y\rest n)\}\le K(n)+c$ for all $n$? (Such a pair is a minimal pair of a particularly simple
kind.) \emph{Negative for c.e.\ reals} \cite[Corollary~1.7]{BarU}; open in general.
\item \cite[Q.~6]{Bar} Is every real reducible to a random real under $\le_{rK}$ or under $\leK$? \emph{The $rK$ part is
answered negatively} \cite{BL}; the $\leK$ part is (L6).
\item \cite[Q.~7]{Bar} Which reals have uncountably many $K$-predecessors? Infinitely often $K$-trivial reals have
countably many \cite{BV}. Barmpalias and Liu \cite{BL2} study the converse in a preprint. They construct a real
that is not infinitely often $K$-trivial and whose $K$-degree (as an equivalence class) is countable, and a minimal
pair of reals that are not infinitely often $K$-trivial \cite[Propositions~1.2--1.3]{BL2}; for suitable order
functions $h$ they construct reals $x$ with uncountably many predecessors such that, writing
$\Delta_K(x,n)=K(x\rest n)-K(n)$, one has $\inf_{m\ge n}\Delta_K(x,m)\le h(n)$ for all large $n$
\cite[Proposition~1.4]{BL2}; and they prove an ample-excess theorem: if $x$ is not infinitely often $K$-trivial,
then $\sum_n2^{\log n-K(x\rest n)}<\infty$ \cite[Theorem~1.5]{BL2}. The cardinality of an equivalence class must be
distinguished from that of the lower cone; whether every real that is not infinitely often $K$-trivial has
uncountably many predecessors is open \cite{BL2}.
\item \cite[Q.~8]{Bar} Is there a maximal $K$-degree (\emph{no}, \cite{Liu})? Is there a minimal $K$-degree
(\emph{open}; Proposition~\ref{prop:dilute} restricts where one could be)?
\item \cite[Q.~9]{Bar} Is there a minimal pair of $\Delta^0_2$ reals in the $K$-degrees? The lowest known complexity is
a $\Sigma^0_2$ real forming a minimal pair with every non-$K$-trivial c.e.\ real \cite{BV}.
\item \cite[Q.~10]{Bar} Are the $K$-degrees (or $rK$-degrees) of c.e.\ sets dense?
\end{enumerate}
Two older questions are settled: there is a minimal pair of $K$-degrees \cite{CM} (asked by Downey and
Hirschfeldt), and there is no minimal pair among the $K$-degrees of c.e.\ reals \cite{BarU} (asked in \cite{DH} and
\cite{MS}).

\subsection{Questions about the structure as a whole}
\begin{enumerate}[label=\textup{(S\arabic*)},leftmargin=3em]
\item Is the first-order theory of $(\mathcal D_K,\le)$ decidable? (Section~\ref{sec:fo}; Remark~\ref{rem:fragments}
for decidable fragments; Propositions~\ref{prop:nies} and~\ref{prop:target} for sufficient conditions for
undecidability.) Is it as complicated as second-order arithmetic, as for the Turing degrees \cite{Simp}?
\item Are the $\forall\exists$ and $\exists\forall$ theories of $(\mathcal D_K,\le)$ decidable? (The latter contains the
sentence ``there is a minimal degree''. For the Turing degrees the $\forall\exists$ theory is decidable and the
$\exists\forall\exists$ theory is not.)
\item Is every finite distributive lattice an initial segment of $\mathcal D_K$ (Question~\ref{q:lat})? Is every
finite lattice a sublattice with least upper bounds and greatest lower bounds preserved
(Question~\ref{q:M3})?
\item Is the class of random degrees first-order definable in $\mathcal D_K$? (Questions~\ref{q1}--\ref{q3}; by
Corollary~\ref{cor:J} the most obvious candidate fails.) Is the class $\{x:K(x\rest n)\le n+O(1)\}$ definable?
\item Does $\mathcal D_K$ have a non-trivial automorphism? Is every automorphism the identity on the random
degrees? (For the Turing degrees the analogous questions are classical and still partly open.)
\item For an oracle $A$, let $\mathcal D_K^A$ be the $K^A$-degrees. Is $\mathcal D_K^A\cong\mathcal D_K$, or at least
$\mathcal D_K^A\equiv\mathcal D_K$, for all $A$? If not, $A\mapsto\mathrm{Th}(\mathcal D_K^A)$ is a degree-invariant
function whose complexity gives lower bounds for $\mathrm{Th}(\mathcal D_K)$.
\item Are the $K$-degrees of left-c.e.\ reals, or of random reals, a definable substructure, and is the theory of
the former decidable? (They form an upper semilattice in which the join is given by addition \cite{FP}.)
\end{enumerate}

\subsection{Can \texorpdfstring{$\mathcal D_K$}{the K-degrees} be isomorphic to another degree structure?}
The answer for the classical structures is no, and for a simple reason.

\begin{proposition}\label{prop:noniso}
\begin{enumerate}[label=\textup{(\alph*)}]
\item $(\mathcal D_K,\le)$ is not elementarily equivalent to any upper semilattice. In particular it is not
elementarily equivalent, let alone isomorphic, to the Turing degrees, the many-one, truth-table, weak
truth-table, enumeration or $LR$-degrees, or to the $K$-degrees of left-c.e.\ reals.
\item $(\mathcal D_K,\le)$ is not isomorphic to any partial order in which every element has countably many
predecessors; this includes the Turing, arithmetic and hyperarithmetic degrees.
\item Every countable partial order embeds into $(\mathcal D_K,\le)$, even below every random degree, and every
finite configuration of $\le$ and $\incomp$ satisfying \textup{(C1)} and \textup{(C2)} occurs. So embeddings of
countable structures do not distinguish $\mathcal D_K$ from the Turing degrees.
\end{enumerate}
\end{proposition}
\begin{proof}
(a) The sentence $\forall x\,\forall y\,\exists z\,(x\le z\wedge y\le z)$ holds in every upper semilattice and fails in
$\mathcal D_K$ by Miller and Yu \cite{MY}. The listed structures are upper semilattices; for the $K$-degrees of
left-c.e.\ reals, where the join is given by addition, see \cite{FP}.
(b) It suffices to find a degree with uncountably many predecessor degrees. For every real $V$ the sparse code
$S_V$ of Section~\ref{sec:sparse} satisfies $K(S_V\rest n)\le n+O(1)$, so $S_V$ is $K$-below every random real
(Proposition~\ref{prop:below}). The lower cone of each $S_V$ is countable (Corollary~\ref{cor:countable}(b)), so the
map $V\mapsto\deg_K(S_V)$ has countable fibres and takes continuum many values. Hence every random degree has
continuum many predecessors.
(c) The Boolean algebra of computable sets modulo finite sets is countable and atomless, and every countable
partial order embeds into it (into the Boolean algebra generated by its principal down-sets, which is countable,
and every countable Boolean algebra embeds into the countable atomless one). Now apply
Theorem~\ref{thm:masks}(a), or Remark~\ref{rem:ideal}, and Corollary~\ref{cor:fin}.
\end{proof}

So the natural candidates for elementary equivalence with $\mathcal D_K$ are the structures that also lack joins:
the $C$-degrees and the vL-degrees of Miller and Yu, in which almost every pair of reals has no upper bound
\cite{MN}. For the vL-degrees the answer is negative by Proposition~\ref{prop:vL}; for the $C$-degrees the
question is open.

\begin{proposition}\label{prop:vL}
Let $\varphi$ be the sentence ``some nonzero degree is compatible with every degree''. Then $\varphi$ holds in
$\mathcal D_K$, assuming the result of \cite{Zhao} used in Proposition~\ref{prop:dim}, and fails in the vL-degrees.
So $\mathcal D_K$ and the vL-degrees are not elementarily equivalent.
\end{proposition}
\begin{proof}
In $\mathcal D_K$: let $R$ be random and let $x$ copy $R$ on the even positions and be $0$ on the odd ones. Then
$K(x\rest n)\le n/2+O(\log n)$, so $x$ is compatible with every real by Proposition~\ref{prop:dim}; and
$K(x\rest n)\ge n/2-O(1)$, because $R\rest\lfloor n/2\rfloor$ is computable from $x\rest n$, so $x$ is not
$K$-trivial.

In the vL-degrees, with $X\le_{vL}Y$ iff for every $Z$, $X\oplus Z$ random implies $Y\oplus Z$ random, all
non-random reals form the least degree, so a nonzero degree is the degree of a random real $X$. Let $Y$ be random relative to $X$. Then $X\oplus Y$ is random, and by Miller and Yu \cite{MY}
$X$ and $Y$ have no upper bound in the vL-degrees. So every nonzero vL-degree is incompatible with some degree.
\end{proof}
Both facts about the vL-degrees used here are proved by Miller and Yu \cite[Theorem~3.4]{MY}.

For restricted structures a separation is known: the $C$-degrees and the $K$-degrees of c.e.\ reals are not
elementarily equivalent, and neither are those of c.e.\ sets \cite[Corollary~1.6]{BarU}. This does not by itself
settle the question for the full structures. For the $C$-degrees, Merkle and Stephan \cite[Theorem~17]{MS} prove
that there is a minimal $C$-degree. So a
negative answer to the minimal $K$-degree question (L15) would show
that $\mathcal D_K$ and the $C$-degrees are not elementarily equivalent; equivalently, elementary equivalence of the
two structures would imply that a minimal $K$-degree exists.

\begin{question}\label{q:C}
Are $(\mathcal D_K,\le)$ and the $C$-degrees elementarily equivalent? Does the sentence $\varphi$ of
Proposition~\ref{prop:vL} hold in the $C$-degrees? Is there any natural degree structure, not defined from
initial segment complexity, that is elementarily equivalent to $\mathcal D_K$?
\end{question}

\subsection{Proposed questions}\label{sec:proposed}
The questions below are, as far as we know, new. Each comes with the reason it would help to settle the decision
problem or to describe $\mathcal D_K$ completely. One preliminary remark fixes the target.

\begin{remark}[An upper bound]\label{rem:upper}
The relation $x\leK y$ on reals is $\Sigma^0_2$, so $\mathcal D_K$ is interpretable in second-order arithmetic: degrees
are reals modulo the $\Sigma^0_2$ relation $\equiv_K$, and a first-order sentence about $\mathcal D_K$ becomes a
sentence of second-order arithmetic with the same quantifier structure over reals. Hence $\mathrm{Th}(\mathcal D_K)$
is computable from the theory of second-order arithmetic. Every question about the complexity of
$\mathrm{Th}(\mathcal D_K)$ is therefore a question about a lower bound.
\end{remark}

\begin{enumerate}[label=\textup{(P\arabic*)},leftmargin=3em]
\item \emph{Joins are maxima?} If $z$ is the least upper bound of $x$ and $y$ in $\mathcal D_K$, is
$K(z\rest n)=\max\{K(x\rest n),K(y\rest n)\}\pm O(1)$? Dually for greatest lower bounds and minima.

\emph{Why.} Lemma~\ref{lem:maxmin} gives the converse. A positive answer would make the partial lattice operations
of $\mathcal D_K$ pointwise operations on functions; then every finite partial lattice in $\mathcal D_K$ whose operations
are defined is distributive, so $N_5$ and $M_3$ do not embed (Question~\ref{q:M3}), and the obstruction to
Route~B in Section~\ref{sec:fo} becomes a question about which pointwise maxima are realized. A negative answer,
for example via Question~\ref{q:lub}, would exhibit joins of a new kind and open room for non-distributive coding.

\item \emph{Realizability of maxima.} For which compatible reals $x,y$ is the function
\[n\mapsto\max\{K(x\rest n),K(y\rest n)\}\]
realized, up to $O(1)$, as the complexity function of a real? Is the set of
such pairs, modulo $\equiv_K$, definable in $\mathcal D_K$?

\emph{Why.} $\mathcal D_K$ is, by definition, the poset of realized complexity functions under pointwise $\le^+$. Its
whole lattice-theoretic behaviour is governed by which pointwise operations stay inside this set; Miller and Yu's
theorem says that some upper bounds are not realized at all.

\item \emph{Density and covers.} If $x<_Ky$, is there $z$ with $x<_Kz<_Ky$? If not, is the set of pairs
$(x,y)$ with $y$ a minimal cover of $x$ definable and large?

\emph{Why.} Density is the first $\forall\exists$ sentence one needs. For the $K$-degrees of left-c.e.\ reals it
holds (Downey, Hirschfeldt and LaForte, as reported in \cite{FP} and \cite{CM}). On sparse codes $\leK$ behaves like
$rK$-reducibility (Corollary~\ref{cor:sparse}), which has a minimal degree \cite{RS}; minimal covers there would
be a source of rigid finite configurations, which is what initial-segment codings need.

\item \emph{The local structure below sparse codes.} For which countable partial orders $L$ is there a real $V$ with
$[0,S_V]\cong L$? In particular, is every finite distributive lattice of this form?

\emph{Why.} These cones are countable, and every element is computable in $V\oplus\emptyset'$
(Proposition~\ref{prop:turing}); so they are local structures like the Turing degrees below $\emptyset'$, where
initial-segment methods are well developed. A positive answer for finite distributive lattices gives
undecidability (Question~\ref{q:lat}). It also forces minimal degrees, so by Proposition~\ref{prop:dilute} the reals $V$
would have to lie in the gap zone of \cite{BaBa}.

\item \emph{Interpreting $rK$.} Is the class $\{\deg_K(S_V):V\in2^\omega\}$ definable in $\mathcal D_K$ with
parameters? Is the theory of the degrees of $V\mapsto S_V$ under the reducibility of Corollary~\ref{cor:sparse}
undecidable?

\emph{Why.} A definable copy of a structure with an undecidable theory gives undecidability. The reducibility of
Corollary~\ref{cor:sparse} is a strong reducibility, close to $rK$, for which tree constructions exist; it may be
far easier to analyse than $\mathcal D_K$ itself.

\item \emph{Exact pairs and definable countable sets.} Is every countable downward-closed set $I$ of degrees that
is directed (any two elements have an upper bound in $I$) of the form $\{x:x\le a\wedge x\le b\}$ for some $a,b$?
Is every countable antichain of degrees below a random degree definable from finitely many parameters?

\emph{Why.} These are the two coding lemmas, Spector's exact pair theorem and the Slaman--Woodin definability of
countable antichains, behind the proof that the Turing degrees interpret second-order arithmetic. Remark~\ref{rem:ideal}
already gives a one-degree version for ideals of masks of computable sets. Positive answers would make
$\mathrm{Th}(\mathcal D_K)$ as complicated as second-order arithmetic, matching Remark~\ref{rem:upper}.

\item \emph{Definability of cardinality.} Is the class of degrees with countably many predecessors definable in
$\mathcal D_K$? Does it coincide with the infinitely often $K$-trivial degrees?

\emph{Why.} Lower cones split into countable ones (where the local structure is arithmetical, by
Proposition~\ref{prop:turing}) and uncountable ones. A definable split would let one study $\mathcal D_K$ as an
arithmetical ``bottom'' plus a large ``top'', and would bear on (S4) and (S5). The second part is the subject of
\cite{BL2}.

\item \emph{The top: upper cones of highly random degrees.} Upper cones of $2$-random degrees are countable
\cite{Mi}, and by \cite{Liu} nonempty. Is the upper cone of a $3$-random degree a chain? Is it well-ordered? Is the
order type the same for almost every real?

\emph{Why.} A complete description of $\mathcal D_K$ must describe its top. If almost every degree has a simple
upper cone, then the random degrees are a structure with an almost-surely constant theory, which may be decidable
even if $\mathcal D_K$ is not.

\item \emph{The theory of the random degrees.} Is the theory of $\mathcal D_K$ restricted to the $1$-random degrees
decidable? Is the random part definable in $\mathcal D_K$ (S4)?

\emph{Why.} Among random reals $\leK$ agrees closely with vL-reducibility \cite{MY}, which is easier to work with.
If the random part is definable and has a decidable theory, all of the complexity of $\mathrm{Th}(\mathcal D_K)$ lies
below the random degrees, and Routes~A--C of this paper are aimed at the right place.

\end{enumerate}

The questions are not independent. A positive answer to (P1) together with (P2) would describe all least upper
bounds; (P4) and (P5) concern the same family of reals; and (P6) together with (P5) or (P4) is a complete
programme for proving that $\mathcal D_K$ is as complicated as second-order arithmetic. In our view (P4) and (P5) are
the best places to start, because they reduce questions about $\mathcal D_K$ to a strong reducibility with
countable cones, where the classical methods of degree theory apply.

\subsection*{Use of AI tools}
The literature review, the proof strategy and first drafts of the proofs and of the text were developed in
interaction with large language models, namely Claude (Anthropic) and models of OpenAI, which were coordinated
through AETHER IDE, a multi-agent environment developed in the author's laboratory
(\url{https://amalithlab.github.io}).


\begin{thebibliography}{99}
\bibitem{AS} K.~Ambos-Spies, \emph{Inhomogeneities in the polynomial-time degrees: the degrees of super sparse sets},
 Inform.\ Process.\ Lett.\ \textbf{22} (1986), 113--117.
\bibitem{BarU} G.~Barmpalias, \emph{Universal computably enumerable sets and initial segment prefix-free
complexity}, Inform.\ and Comput.\ \textbf{233} (2013), 41--59.
\bibitem{Bau} B.~Bauwens, \emph{Relating and contrasting plain and prefix Kolmogorov complexity},
arXiv:1311.2092, 2013.
\bibitem{BaBa} M.~Baartse and G.~Barmpalias, \emph{On the gap between trivial and nontrivial initial segment
prefix-free complexity}, Theory Comput.\ Syst.\ \textbf{52} (2013), 28--47.
\bibitem{Bar} G.~Barmpalias, \emph{Algorithmic randomness and measures of complexity},
 Bull.\ Symbolic Logic \textbf{19} (2013), 318--350.
\bibitem{BL} G.~Barmpalias and L.~Liu, \emph{Martin-L\"of randomness and $rK$-reducibility}, arXiv:2610.39174, 2026, preprint.

\bibitem{BL2} G.~Barmpalias and L.~Liu, \emph{Cardinality of the $K$-cones below a non infinitely often $K$-trivial},
 arXiv:2610.39172, 2026, preprint.

\bibitem{BV} G.~Barmpalias and C.~S.~Vlek, \emph{Kolmogorov complexity of initial segments of sequences and
 arithmetical definability}, Theoret.\ Comput.\ Sci.\ \textbf{412} (2011).
\bibitem{Blog17} A.~Nies (ed.), \emph{Logic Blog 2017}, arXiv:1804.05331.
\bibitem{CM} B.~F.~Csima and A.~Montalb\'an, \emph{A minimal pair of $K$-degrees}, Proc.\ Amer.\ Math.\ Soc.\
\textbf{134} (2006), 1499--1502.
\bibitem{DH} R.~G.~Downey and D.~R.~Hirschfeldt, \emph{Algorithmic Randomness and Complexity}, Springer, 2010.
\bibitem{DN} R.~Downey and A.~Nies, \emph{Undecidability results for low complexity time classes},
 J.~Comput.\ System Sci.\ \textbf{60} (2000), no.~2.
\bibitem{DHL} R.~Downey, D.~R.~Hirschfeldt and G.~LaForte, \emph{Undecidability of the structure of the Solovay degrees
 of c.e.\ reals}, J.~Comput.\ System Sci.\ \textbf{73} (2007), 769--787.
\bibitem{ELTT} Yu.~L.~Ershov, I.~A.~Lavrov, A.~D.~Taimanov and M.~A.~Taitslin, \emph{Elementary theories},
 Russian Math.\ Surveys \textbf{20} (1965), no.~4, 35--105.
\bibitem{FP} J.~N.~Y.~Franklin and C.~P.~Porter, \emph{Key developments in algorithmic randomness},
arXiv:2004.02851; in: Algorithmic Randomness: Progress and Prospects, Lecture Notes in Logic 50, 2020.
\bibitem{Lach} A.~H.~Lachlan, \emph{Distributive initial segments of the degrees of unsolvability},
 Z.~Math.\ Logik Grundlagen Math.\ \textbf{14} (1968), 457--472.
\bibitem{Liu} L.~Liu, \emph{There is no maximal $K$-degree}, arXiv:2609.11054, 2026, preprint.

\bibitem{ML} P.~Martin-L\"of, \emph{Complexity oscillations in infinite binary sequences},
 Z.~Wahrscheinlichkeitstheorie verw.\ Gebiete \textbf{19} (1971), 225--230.
\bibitem{MS} W.~Merkle and F.~Stephan, \emph{On $C$-degrees, $H$-degrees and $T$-degrees}, in: Proc.\ 22nd IEEE
 Conference on Computational Complexity (CCC 2007), 60--69.
\bibitem{Mi} J.~S.~Miller, \emph{The $K$-degrees, low for $K$ degrees, and weakly low for $K$ sets},
Notre Dame J.~Formal Logic \textbf{50} (2009), 381--391.
\bibitem{MN} J.~S.~Miller and A.~Nies, \emph{Randomness and computability: open questions},
 Bull.\ Symbolic Logic \textbf{12} (2006), 390--410.
\bibitem{MY} J.~S.~Miller and L.~Yu, \emph{On initial segment complexity and degrees of randomness},
 Trans.\ Amer.\ Math.\ Soc.\ \textbf{360} (2008), 3193--3210.
\bibitem{Mu} An.~A.~Muchnik, \emph{Conditional complexity and codes}, Theoret.\ Comput.\ Sci.\ \textbf{271}
(2002), 97--109.
\bibitem{MRS} D.~Musatov, A.~Romashchenko and A.~Shen, \emph{Variations on Muchnik's conditional complexity
theorem}, Theory Comput.\ Syst.\ \textbf{49} (2011), 227--245.
\bibitem{Nies97} A.~Nies, \emph{Intervals of the lattice of computably enumerable sets and effective Boolean algebras},
 Bull.\ London Math.\ Soc.\ \textbf{29} (1997), 683--692.
\bibitem{Nies00} A.~Nies, \emph{Effectively dense Boolean algebras and their applications},
 Trans.\ Amer.\ Math.\ Soc.\ \textbf{352} (2000), 4989--5012.
\bibitem{Rabin} M.~O.~Rabin, \emph{Decidability of second-order theories and automata on infinite trees},
 Trans.\ Amer.\ Math.\ Soc.\ \textbf{141} (1969), 1--35.
\bibitem{RS} A.~Raichev and F.~Stephan, \emph{A minimal $rK$-degree}, in: Computational Prospects of Infinity, Part~II,
 World Scientific, 2008.
\bibitem{Silver} J.~Silver, \emph{Every analytic set is Ramsey}, J.~Symbolic Logic \textbf{35} (1970), 60--64.
\bibitem{Simp} S.~G.~Simpson, \emph{First-order theory of the degrees of recursive unsolvability},
 Ann.\ of Math.\ \textbf{105} (1977), 121--139.
\bibitem{Zhao} X.~Zhao, \emph{Effective Hausdorff dimension one is cofinal in the $K$-degrees}, arXiv:2611.01421, 2026, preprint.
\end{thebibliography}
\end{document}